\documentclass[preprint]{elsarticle}

\usepackage{style}
\usepackage{macros}

\usepackage{mathptmx}

\SetMathAlphabet{\mathcal}{normal}{OMS}{lmsy}{m}{n}
\SetMathAlphabet{\mathcal}{bold}{OMS}{lmsy}{m}{n}

\usepackage{autonum}

\usepackage[left=1in,right=1in]{geometry}

\usepackage[bf,small,raggedright]{titlesec}

\title{H\"older estimates for fractional parabolic equations\\ with critical divergence free drifts} 
\author{Mat\'ias G. Delgadino}
\address{ICTP, Strada Costiera 11, Trieste, Italy, 34151}
\author{Scott Smith}
\address{MIS-MPG, Inselstra{\ss}e 22, Leipzig,
Germany, 04103}
\begin{document}
\begin{abstract}
This work focuses on drift-diffusion equations with fractional dissipation $(-\Delta)^{\alpha}$ in the regime $\alpha \in (1/2,1)$.  Our main result is an a priori H\"older estimate on smooth solutions to the Cauchy problem, starting from initial data with finite energy.  We prove that for some $\beta \in (0,1)$, the $C^{\beta}$ norm of the solution depends only on the size of the drift in critical spaces of the form $L^{q}_{t}(BMO^{-\gamma}_{x})$ with $q>2$ and $\gamma \in (0,2\alpha-1]$, along with the $L^{2}_{x}$ norm of the initial datum.  The proof uses the Caffarelli/Vasseur variant of De Giorgi's method for non-local equations. 
\end{abstract}
\maketitle

\tableofcontents
\section{Introduction}\label{S.Intro}
This paper is concerned with quantitative estimates for solutions to the following partial differential equation:
\begin{equation}\label{eq:DriftDiff}
\partial_{t}\theta + u \cdot \nabla \theta +(-\Delta)^{\alpha}\theta = 0 \quad \text{in} \quad \R_{+}\times \R^{d}.
\end{equation}
The quantity $\theta(t,x)$ is a scalar and the $\R^{d}$ valued vector field $u(t,x)$ is time dependent and divergence free.  We study the Cauchy problem, and supplement \eqref{eq:DriftDiff} with initial datum $\theta_{0}$ in the natural energy space $L^{2}(\R^{d})$.  Under the qualitative hypothesis that $u$ and $\theta_{0}$ are smooth, \eqref{eq:DriftDiff} admits a unique classical solution $\theta$ starting from $\theta_{0}$.  We are interested in precisely quantifying the regularity of the solution in terms of the drift and the initial data.  More specifically, our aim is to obtain an a priori estimate for the H\"older norm of $\theta$ in terms of the weakest affordable norm of $u$, together with the initial energy.

The modern approach to this question begins with a classification of drifts according to their criticality.  At the level of a qualitative heuristic, criticality refers to the relative strength of the advection $u \cdot \nabla \theta $ versus the diffusion $(-\Delta)^{\alpha} \theta$ on the small scales.  In the literature on active scalars,  $u$ often has a fixed functional relationship with $\theta$, and the criticality of equation \eqref{eq:DriftDiff} changes by varying $\alpha$.  In contrast, the perspective here is that $\alpha$ is fixed, and criticality varies with the roughness of the drift $u$.  For sub-critical drifts, diffusion dominates and one expects \eqref{eq:DriftDiff} to obey continuity estimates similar to the fractional heat equation.  However, for super-critical drifts advection dominates and \eqref{eq:DriftDiff} can behave more like a transport equation, where initial discontinuities propagate in time.  The last alternative is that $u$ is critical, meaning that both influences are equally balanced on small scales.  This case requires careful analysis, and our work focuses primarily on this regime.

To quantify the heuristics given above, we consider the interplay between the natural scaling transformation preserving \eqref{eq:DriftDiff} and the norm measuring the drift.  Given a solution $\theta$ and a drift $u$, the re-scaling
\begin{equation}\label{eq:Scaling_Transformation}
\theta_{\lambda}(t,x)=\theta(\lambda^{2\alpha}t,\lambda x), \quad u_{\lambda}(t,x)=\lambda^{2\alpha-1}u(\lambda^{2\alpha}t,\lambda x)
\end{equation}
produces a new solution $\theta_{\lambda}$ relative to the drift $u_{\lambda}$ and the scale $\lambda>0$.  Criticality can now be quantified in terms of linear spaces $X$ endowed with a homogeneous norm.  Namely, $X$ is said to be sub-critical (or super-critical) if $\|u_{\lambda}\|_{X}$ tends to zero (or infinity) as $\lambda \to 0$.  In contrast, the space is called critical if $\|u_{\lambda}\|_{X}$ is independent of $\lambda$.

The particular critical spaces studied in this article are based on negative $\text{BMO}(\R^{d})$ norms in space, which we now describe.  For each time $t>0$, $u(t)$ is taken in a space denoted $\text{BMO}^{-\gamma}(\R^{d};\R^{d})$ with $\gamma \in (0,1)$. One can think of this as a vector space of distributions with negative order, obtained by ``taking $\gamma$ derivatives'' of a vector field in $\text{BMO}$. We defer to Section \ref{Sec:NegBMO} for a precise definition, only noting this hypothesis yields a $\psi(t) \in \text{BMO}(\R^{d};\R^{d})$ such that $u(t)=(-\Delta)^{\frac{\gamma}{2}}\psi(t)$ componentwise.  Moreover, the corresponding norm is defined by
\begin{equation} \label{eq:Intro_Neg_BMO_Norm}
\|u(t)\|_{\text{BMO}^{-\gamma}(\R^{d};\R^{d}) }=\|\psi(t)\|_{\text{BMO}(\R^{d};\R^{d})}.
\end{equation}  
In fact, we impose that $t \to \|u(t)\|_{\text{BMO}^{-\gamma}(\R^{d};\R^{d})}$ belongs to $L^{q}(\R_{+})$ and use the shorthand notation $u \in \LB{q}{-\gamma}$.  This defines a critical space for the drift provided the exponents satisfy
\begin{equation}\label{eq:Intro_Crit_Exponents}
\frac{2\alpha}{q}+\gamma=2\alpha -1, \quad \quad \gamma \in (0,2\alpha-1].
\end{equation}
Next, we introduce a set $\mathcal{S}(D,\alpha)$ consisting of all possible smooth solutions to \eqref{eq:DriftDiff} with critical drift norm at most $D>0$, relative to the order of dissipation $\alpha \in (1/2,1)$.  More precisely, $\theta \in \mathcal{S}(D,\alpha)$ provided there exists a smooth, divergence free $u \in\LB{q}{-\gamma}$ such that \eqref{eq:DriftDiff} holds classically and
\begin{equation}
\|u\|_{\LB{q}{-\gamma}}\leq D,
\end{equation}
for some exponents $\gamma, \, q$ satisfying \eqref{eq:Intro_Crit_Exponents}.  
Finally, we define a H\"older semi-norm which is consistent with the scaling \eqref{eq:Scaling_Transformation} via
\begin{equation}\label{eq:HolderSemiNorm}
[\theta]_{C^{\beta}_{\alpha}([t,\infty) \times \R^{d})}=\sup_{(s,x),(r,y) \in [t,\infty] \times \R^{d}}\frac{|\theta(s,x)-\theta(r,y)|}{|s-r|^{\frac{\beta}{2\alpha}}+|x-y|^{\beta}}.
\end{equation}
Our main result is the following a priori bound for the H\"older semi-norm of $\theta$ in $C^{\beta}_{\alpha}$.
\begin{Thm} \label{Thm:HolderContResult}
For all $\alpha \in (1/2,1)$ and $D>0$, there exist positive constants $\beta, C$ such that for all $\theta \in \mathcal{S}(D,\alpha)$ and times $t>0$ , 
\begin{equation}
[\theta]_{C^{\beta}_{\alpha} ( [t,\infty) \times \R^{d} )} \leq C t^{-(\frac{d}{4\alpha}+\beta) }\|\theta_{0}\|_{L^{2}(\R^{d})}.
\end{equation}
Moreover, the constants $\beta$ and $C$ are universal (depend only on $D$,$\alpha$, and $d$).
\end{Thm}

Before comparing Theorem \ref{Thm:HolderContResult} to the existing literature, let us return to our discussion of critical spaces, focusing in particular on those of the form $L^q_t(Y_x)$ for $q \in [1,\infty]$ and $Y_{x} \subset \mathcal{D}'(\R^{d};\R^{d})$.  A first remark is that these spaces vary substantially as the level of dissipation ranges from the transport regime $\alpha=0$ to the case of full diffusion $\alpha=1$. As $\alpha$ grows, the increased diffusion can be used to compensate for more irregular drifts and the critical spaces become larger.  In particular, an important transition occurs at $\alpha=1/2$.  Namely, for $\alpha \in (0,1/2]$, $Y_{x}$  can only include distributions with non-negative order, but for $\alpha \in (1/2,1]$, $Y_{x}$ can include drifts with negative regularity.  The works closest to ours concern the divergence free context, so we describe these first. The question of removing this assumption is postponed until the end of the introduction. The literature in the critical case is devoted to one of two scenarios: the first is the regime $\alpha \in (0,1/2]$ and the second $\alpha=1$.  

In the case $\alpha \in (0,1/2]$, the seminal work is by Caffarelli-Vasseur \cite{caffarelli2010drift}, which employs a De Giorgi scheme in the case $\alpha=1/2$ to treat drifts in $L^{\infty}_{t}(\text{BMO}_{x})$.  Subsequently,  Constantin-Wu \cite{constantin2008regularity} prove that the method of Caffarelli-Vasseur can also be applied in the regime $\alpha \in (0,1/2)$, provided $u$ belongs to $L^{\infty}_{t}(C^{1-2\alpha}_{x})$.  In these two articles, the main obstacle is the lack of a true local energy inequality, due to the presence of the fractional Laplacian.  To address this problem, the authors of \cite{caffarelli2010drift} developed a technique for embedding the non-local obstructions into the solution of an auxiliary problem, the harmonic extension of the solution to the upper half plane.  They prove a form of the Caccioppoli inequality involving both the solution $\theta$ and its harmonic extension $\theta^{*}$, then successively exploit the estimate in a sophisticated De Giorgi scheme, yielding a decay of the oscillations. In this article, we rely heavily on the methods developed in \cite{caffarelli2010drift}.

In the case $\alpha=1$, Osada \cite{osada1987diffusion} treated the critical case $u \in L^{\infty}_{t}(W^{-1,\infty}_{x})$ nearly thirty years ago.  More recently, two different sets of authors studied drifts in the more refined space $\LB{\infty}{-1}$, each using different methods.  In \citep{seregin2012divergence}, Seregin-Silvestre-\u{S}ver\'{a}k-Zlato\u{s} establish the desired H\"older bound by means of a Harnack inequality proven with Moser's iteration technique.
In \cite{friedlander2011global}, Friedlander-Vicol work more in the spirit of \cite{caffarelli2010drift}, developing a De Giorgi scheme and proving an oscillation reduction result.  The main difficulty in these articles is dealing with the negative regularity of the drift.  A key problem the authors overcame in \cite{friedlander2011global} was obtaining a form of Caccioppoli inequality with constants depending only on the $\LB{\infty}{-1}$ norm of the velocity.  In fact, they obtain an estimate which is weaker than the Caccioppoli inequality for the heat equation, but show nonetheless that it is sufficent for the purpose of the De Giorgi scheme.

It is important to emphasize that the key obstacles faced in the regime $\alpha \in (0,1/2]$ are somewhat distinct from the case $\alpha=1$.  While the first deals mostly with the non-locality, the second deals primarily with the negative regularity of the drift.  In the case $\alpha \in (1/2,1)$ we face both these difficulties simultaneously.  

A key part of our work is to obtain a form of Caccioppoli inequality involving only the rough norms of the drift.  This is the content of Sections \ref{S.ProdEst} and \ref{S.Cacciop}.  The proof is inspired by a trick in \citep{friedlander2011global}, where an integration by parts combined with an application of the John-Nirenberg inequality allows to locally transfer the negative regularity of the drift onto the solution, then use the dissipative bounds.  In the present setting, we cannot afford to pass off a full derivative onto $\theta$, so we use a fractional integration by parts instead.  This is slightly subtle because it turns a localized integral into an integral over the whole space.  Nonetheless, in Lemma \ref{L.BMOProd} we obtain our basic mechanism for trading fractional negative regularity off of the drift and onto the solution at a global level.  

However, for the Caccioppoli estimate we must perform this locally.  As in \cite{caffarelli2010drift} and \cite{constantin2008regularity}, our substitute for the local energy inequality uses a solution to an auxilliary problem, given by the Caffarelli-Silvestre \cite{caffarelli2007extension} extension.  In Lemma \ref{L.DriftEstimate}, we show how to control the drift contribution in terms of both the extension and the true solution.  Due to the lack of a pointwise fractional product rule, this estimate costs a small error term, which we must iterate away in Section \ref{S.Cacciop}. Ultimately this leads to a slightly less precise form of Caccioppoli estimate than in the work of Caffarelli-Vasseur. However, in some sense it is stronger than the estimate obtained in \citep{friedlander2011global}, and we explicitly use the strict inequality $\alpha<1$.  

The remainder of the paper is devoted to applying the Caccioppoli estimate together with the method of Caffarelli-Vasseur to obtain the H\"older estimate through a De Giorgi iteration. The reader is referred to \citep{caffarelli2010giorgi} for an introduction to this technique.  Let us finish with some comments on the sharpness of the hypotheses of Theorem \ref{Thm:HolderContResult}.

The results of \cite{silvestre2013loss} show that one cannot generally expect Theorem \ref{Thm:HolderContResult} to hold if the drift norm is only controlled in a super-critical space.  In particular, this means that if we focus on critical norms of the form $\LB{q}{-\gamma}$, then the exponents in \eqref{eq:Intro_Crit_Exponents} are optimal.  Indeed, even though the results in \cite{silvestre2013loss} are stated for vector fields in $L^{\frac{d}{2\alpha-1}}(\R^{d};\R^{d})$, this space embeds into $\text{BMO}^{1-2\alpha}(\R^{d};\R^{d})$.  Thus, the results in \cite{silvestre2013loss} prevent Theorem \ref{Thm:HolderContResult} due to the endpoint case $q=\infty$.  As stated, the results in \cite{silvestre2013loss} do not rule out Theorem \ref{Thm:HolderContResult} if we consider only $q \in (2,\infty)$, but we are inclined to believe there are counterexamples in this regime also.  On the other hand, it is conceivable that there are larger critical spaces where Theorem \ref{Thm:HolderContResult} could hold.  The next natural question is whether Theorem~\ref{Thm:HolderContResult} can be extended to the Besov space $L^{q}_{t}(B^{-\gamma}_{\infty,\infty})$.

Let us also comment on the the divergence free hypothesis for the velocity. The work of Silvestre \cite{silvestre2010holder,silvestre2010differentiability} shows that this assumption can be dropped if $\alpha \in (0,1/2)$, provided one forgoes the $\text{BMO}$ framework, working in a bit larger critical spaces.  Namely, for $\alpha \in (0,1/2)$ the result of \cite{constantin2008regularity} is shown to hold without the divergence free assumption.  In the case $\alpha \in [1/2,1)$, H\"older estimates are obtained if the drift is bounded.  In particular, when $\alpha=1/2$ the results of Caffarelli-Vasseur can almost be recovered without the divergence free assumption, but not quite up to drifts in $\text{BMO}$.    It is important to note that \cite{silvestre2010holder,silvestre2010differentiability} do not follow a variational approach.  

The techniques in this paper seem to break down if one removes entirely the divergence free hypothesis.  However, it is likely that the one could at least extend our results to the case $\Div u \in \LB{\infty}{-\alpha}$.  Note that belonging to the critical space $\LB{\infty}{1-2\alpha}$ only imposes $\Div u \in \LB{\infty}{-2\alpha}$, so there is a significant gap that remains.

Finally, we want to mention the possibility of building weak solutions (in the variational sense) to \eqref{eq:DriftDiff} when $u$ is truly a distribution in $\LB{q}{-\gamma}$.  If $u$ is divergence free, or more generally $\Div u \in \LH{2}{-\alpha}$, then one can define $u \cdot \nabla \theta$ as a distribution by way of a formal fractional integration by parts.  One could construct H\"older continuous weak solutions of this type using Theorem \ref{Thm:HolderContResult}.  We will discuss these matters in a separate article, where results on existence and uniqueness for fractional parabolic equations with rough coefficients will be proved in the spirit of Di-Perna-Lions \cite{diperna1989cauchy} and Lions-Le-Bris \cite{bris2008existence, le2004renormalized}.  In particular, we intend to discuss the relation between criticality and criteria for uniqueness.

\section{Preliminaries}\label{S.Prelim}
\subsection{Negative BMO spaces}\label{Sec:NegBMO}
Let us begin by introducing the negative BMO space $\BMO{-\gamma}$.
\begin{Def}
For $\gamma>0$, a distributional vector field $u \in \mathcal{D}'(\R^{d};\R^{d})$ belongs to $\BMO{-\gamma}$ provided there exists $\psi \in \text{BMO}(\R^{d};\R^{d})$ such that $u =(-\Delta)^{\frac{\gamma}{2}} \psi$ componentwise.
\end{Def}
Observe that for $\psi \in \text{BMO}(\R^{d};\R^{d})$, the tempered distribution $(-\Delta)^{\frac{\gamma}{2}} \psi \in \mathcal{S}'(\R^{d})$ is defined componentwise by duality for $\varphi \in \mathcal{S}(\R^{d})$ via
\begin{equation}\label{eq:BMO_Duality_Def}
 \left\langle (-\Delta)^{\frac{\gamma}{2}}\psi_{i} ,  \varphi \right\rangle = \int_{\R^{d}}\psi_{i}(x)(-\Delta)^{\frac{\gamma}{2}}\varphi(x)\dx.
\end{equation}
The integral in \eqref{eq:BMO_Duality_Def} is well defined in view of the following two classical inequalities
\begin{equation}
\sup_{x \in \R^{d}} \left[\left(1+|x|^{d+\gamma} \right)|(-\Delta)^{\frac{\gamma}{2}}\varphi(x)| \right ] \leq C, \quad \quad \int_{\R^{d}}\frac{|\psi(x)|}{1+|x|^{d+\gamma}}\dx < \infty.
\end{equation} 
The first is a classical fact about the behavior of $(-\Delta)^{\frac{\gamma}{2}}$ on the Schwartz space $\mathcal{S}(\R^{d})$, while the second is property of functions in $\text{BMO}(\R^{d})$ (see \cite{stein2016harmonic}).  Next we note that $\BMO{-\gamma}$ is a linear space and a norm can be defined as follows: 
\begin{equation}\label{eq:BMO-gammanorm}
\|u\|_{\BMO{-\gamma}}=\|\psi\|_{\text{BMO}(\R^{d};\R^{d})},
\end{equation}
where $\psi$ satisfies $(-\Delta)^{\frac{\gamma}{2}}\psi = u$.  The following Lemma ensures the norm is well-defined.
\begin{Lem}\label{lem:uniqueBMO}
For each $u \in \mathcal{S}'(\R^{d};\R^{d})$, there is at most one solution $\psi \in \text{BMO}(\R^{d};\R^{d})$ to the equation
\begin{equation}
(-\Delta)^{\frac{\gamma}{2}}\psi=u \quad \text{in} \quad \mathcal{D}'(\R^{d};\R^{d}),
\end{equation}
modulo shifts by a constant vector.
\end{Lem}

\begin{proof}
Let $\psi_{1},\psi_{2} \in \text{BMO}(\R^{d};\R^{d})$ satisfy $(-\Delta)^{\frac{\gamma}{2}} \psi_{1}= (-\Delta)^{\frac{\gamma}{2}} \psi_{2} $ in $\mathcal{D}'(\R^{d};\R^{d})$.  Defining the difference $q=\psi_{1}-\psi_{2}$, it follows that each component $q_{i}$ is an $\gamma$ harmonic function satisfying
\begin{equation} \label{eq:BMOGrowth}
\int_{\R^{d}}\frac{|q_{i}(x)|}{1+|x|^{d+\gamma}}\dx < \infty.
\end{equation}
Combining this with the Liouville Theorem for $(-\Delta)^{\frac{\gamma}{2}}$ proved in \cite[Theorem 1.3]{chen2015some}, it follows that each $q_{i}$ is a constant function. Hence, $\psi_{1}$ differs from $\psi_{2}$ by a constant vector.
\end{proof}
\begin{Rem}
Observe that \eqref{eq:BMO-gammanorm} defines a norm rather than a semi-norm. This follows from the observation that if $\|u\|_{\text{BMO}^{-\gamma}(\R^{d})}=0$, then there exist a constant function $\psi$ such that $u = (-\Delta)^{\frac{\gamma}{2}}\psi$, therefore $u=0$.
\end{Rem}
In addition, we define the space of vector valued functions $\LB{q}{-\gamma}$ as follows:
\begin{Def}
An element $u \in \LB{q}{-\gamma}$ is an equivalence class of measurable maps $u: \R_{+} \to \mathcal{D}'(\R^{d};\R^{d})$ such that $t \to \|u(t)\|_{\BMO{-\gamma}}$ belongs to $L^{q}(\R_{+})$.
\end{Def}

\subsection{The solution set}
For each $\alpha \in (1/2,1)$ and $D>0$, we define a set $\mathcal{S}(D,\alpha)$ of solutions to \eqref{eq:DriftDiff}.  The parameter $D$ measures the size of the drift in a critical $\LB{q}{-\gamma}$ space.  Since we only aim in this paper to establish a priori bounds, we restrict attention to classical solutions $\theta$ to $\eqref{eq:DriftDiff}$ driven by regular velocity fields.  By classical, we mean that $\eqref{eq:DriftDiff}$ holds pointwise for $(t,x) \in \R_{+} \times \R^{d}$.  The precise definition of $\mathcal{S}(D,\alpha)$ is as follows:
\begin{Def}
A function $\theta:[0,\infty) \times \R^{d} \to \R$ belongs to $\mathcal{S}(D,\alpha)$ provided it satisfies \eqref{eq:DriftDiff} classically, relative to a smooth, divergence free velocity field $u \in \LB{q}{-\gamma}$ with 
\begin{equation}
\|u\|_{\LB{q}{-\gamma}} \leq D,
\end{equation}
where $\gamma >0$ and $q \geq 1$ are related by $\frac{2\alpha}{q}+\gamma=2\alpha -1$. 
\end{Def}
Given $\theta \in \mathcal{S}(D,\alpha)$ and $\lambda>0$, defining $\theta_{\lambda}$ via \eqref{eq:Scaling_Transformation} yields another element of $\mathcal{S}(D,\alpha)$, relative to the velocity field $u_{\lambda}$.  This is a consequence of the relation \eqref{eq:Intro_Crit_Exponents} defining the exponents $q,\gamma$.

Observe that we do not place any restrictions on the initial value of solutions $\theta \in\mathcal{S}(D,\alpha)$.  In particular, $\theta \in \mathcal{S}(D,\alpha)$ does not imply $\theta(0) \in L^{2}(\R^{d})$.  This ensures that for $\theta \in \mathcal{S}(D,\alpha)$ and a constant $c \in \R$, we have $\theta + c \in \mathcal{S}(D,\alpha)$.

\subsection{The Caffarelli-Silvestre extension problem}
Here we introduce the Caffarelli-Silvestre extension that appears systematically in our work. The reader should be advised that, to avoid technical issues, the results presented in this section are tailored to the needs of the present work and are not stated in full generality. The interested reader is referred to \citep{caffarelli2007extension} and \citep{capella2011regularity}.

We begin by defining a natural weighted norm for the extension, based on the weight $\omega_{\alpha}(z)=z^{1-2\alpha}$ defined for $z \in \R_{+}$.  The space defined below appears repeatedly throughout the article.
\begin{Def}\label{def:weightedspaces}
Given $f\in L^1_{loc}(\R^d\times \R_+)$, we define the $L^2_{\omega_\alpha}(\R^d\times \R_+)$ norm by
\begin{equation}
\|f \|_{L^2_{\omega_\alpha}(\R^d\times \R_+)}^2
=\int_0^\infty\int_{\R^d} z^{1-2\alpha}|f|^2 \dee x \dee z.
\end{equation}
\end{Def}

Next, we define the Caffarelli-Silvestre extension (see \cite{caffarelli2007extension}).
\begin{Def}\label{def:defof*}
Given $f\in H^\alpha(\R^d)$, we define for $(x,z) \in \R^d\times [0,\infty)$,
\begin{equation}\label{eq:definitionof*}
f^{*}(x,z)=f* P_{z}^{\alpha}(x),
\end{equation}
where
\begin{equation}\label{eq:poissonkernel}
P_{z}^{\alpha}(x)=c_{d,\alpha}\frac{z^{2\alpha}}{(|x|^{2}+z^{2})^{\frac{d}{2}+\alpha}}=\frac{1}{z^{d}}P_{1}^{\alpha}\left(\frac{x}{z}\right),
\end{equation}
and $c_{d,\alpha}$ is a normalizing constant such that $\|P_{1}^{\alpha}\|_{L^1(\R^d)}=1$.
\end{Def}
\begin{Rem}
As $\alpha$ is fixed throughout the paper, there is no ambiguity in the definition of $*$.
\end{Rem}
\begin{Rem}\label{rem:uniqueness}
Formally, $f^*$ solves
\begin{equation}\label{eq:extension}
\begin{cases}
\Div(z^{1-2\alpha}\nabla f^{*})=0 &\quad \text{in} \quad \R^{d} \times \R_{+} \\
f^*(x,0)=f(x) &\quad \text{in} \quad \R^{d}, \\
\end{cases}
\end{equation}

In fact, the solution to \eqref{eq:extension} is unique. This can be seen from the fact that \eqref{eq:extension} is also the Euler-Lagrange condition associated to the energy
\begin{equation}\label{eq:energy}
\frac{1}{2}\|\nabla g\|^2_{L^2_{\omega_\alpha}(\R^d\times\R_+)}=\int_0^\infty\int_{\R^d} z^{1-2\alpha}|\nabla g(x,z)|^2\;\dx\dz.
\end{equation}
Uniqueness follows from the fact that \eqref{eq:energy} is, up to constants, strictly convex.
\end{Rem}




A fundamental tool for our analysis is the following key identity, lifting a local average of the fractional laplacian to a weighted inner product of local derivatives on a space of dimension one degree higher. 
\begin{Lem} \label{L.FundamentalIdentity}
For all $f\in H^{\alpha}(\R^{d})\cap W^{1,\infty}(\R^{d})$ and $g\in W^{1,\infty}_c(\R^d\times[0,\infty))$ the following identity holds 
\begin{equation}\label{eq:keyidentity}
\int_{\R^{d}}(-\Delta)^{\alpha}f(x)g_{0}(x)\dx = \int_{\R^{d}}\int_{0}^{\infty}z^{1-2\alpha}\nabla f^{*} \cdot \nabla g\; \dx\dz, 
\end{equation}
where $g_0=g(\cdot,0)$ and $f^*$ is the Caffarelli-Silvestre extension of $f$.
\end{Lem}
\begin{proof}
We introduce two functionals $G$ and $H$ defined by
\begin{align}\label{eq:2ndenergy}
&G(\rho)=\frac{1}{2}\|\nabla \rho\|_{L^2_{\omega_\alpha}(\R^d\times\R_+)}-\int_{\R^{d}}\rho_0(-\Delta)^{\alpha}f \dx. \\
&H(\rho)=\frac{1}{2}\|\rho_0\|_{H^\alpha(\R^d)}^2-\int_{\R^d}\rho_0(-\Delta)^{\alpha}f \dx.
\end{align}
Both functionals are stictly convex over the the convex constraint set
\begin{equation}
\mathcal{A}=\left\{\rho: \R^{d} \times \R_{+} \to \R \mid \nabla \rho \in L^{2}_{\omega_{\alpha}}(\R^{d} \times \R_{+}), \quad \rho_{0} \in H^{\alpha}(\R^{d}) \right\},
\end{equation}
where $\rho_{0}(x)= \rho(x,0)$ for $x \in \R^{d}$.  

The main observation is that \eqref{eq:keyidentity} is the Euler-Lagrange condition for the unique minimizer of $G$ in $\mathcal{A}$, which we want to show is $f^{*}$.  Taking perturbations supported on $\R^d\times (0,\infty)$, the first order conditions imply that any minimizer must solve
\begin{equation}
\nabla\cdot(z^{1-2\alpha}\nabla \rho)=0\qquad\mbox{in} \quad \R^d\times\R_+.
\end{equation}
Since the minimizer belongs to $\mathcal{A}$, it must be the unique Caffarelli-Silvestre extension of it's trace.  Moreover, if $\rho \in \mathcal{A}$ is given by $\rho_{0}^{*}$, then $G(\rho)=H(\rho)$ in view of the identity
\begin{equation}
\|\nabla (\rho_{0}^{*})\|_{L^2_{\omega_\alpha}(\R^d\times\R_+)}=\|\rho_{0}\|_{H^{\alpha}(\R^{d})},
\end{equation}
which has been shown in Section 3.2 of \cite{caffarelli2007extension}.  Hence the minimizer of $G$ in $\mathcal{A}$ also minimizes $H$ in $\mathcal{A}$.  Moreover, $H$ can be re-written as
\begin{equation}
H(\rho)=\frac{1}{2}\|\rho_{0}-f\|_{H^{\alpha}(\R^{d})}-\frac{1}{2}\|f\|_{H^{\alpha}(\R^{d})}.
\end{equation}
Hence, the minimizer of $H$ over $\mathcal{A}$ is clearly $\rho=f^{*}$, completing the proof.
\end{proof}
Using the previous Lemma, we can give a variational characterization to the $H^\alpha(\R^d)$ seminorm in the extended space.
\begin{Cor}\label{prop:variational}
Given $\alpha\in (0,1)$ and $f\in W^{1,\infty}(\R^d) \cap H^{\alpha}(\R^{d})$, 
\begin{equation}\label{eq:VariationalHalpha}
\int_{\R^d}|(-\Delta)^\frac{\alpha}{2}f|^2\;\dx =[f]_{\dot{H}^\alpha(\R^d)}^2=\inf_{g\in S(f)}\int_0^\infty\int_{\R^d} z^{1-2\alpha}|\nabla g(x,z)|^2\;\dx\dz,
\end{equation}
where $\mathcal{A}(f)=\{g\in W^{1,\infty}(\R^d\times[0,\infty))\;: \, \nabla g \in L^{2}_{\omega_{\alpha}}(\R^{d} \times \R_{+}), \, \, \, \;g(x,0)=f\}$.

Moreover, the infimum in \eqref{eq:VariationalHalpha} is attained by $f^*$.
\end{Cor}

\section{Product Estimates}\label{S.ProdEst}
For divergence free drifts $u$, the global balance of energy for \eqref{eq:DriftDiff} is
\begin{equation}\label{eq:GlobalEnergyIneq}
\frac{d}{dt}\|\theta(t)\|_{L^{2}(\R^{d})}^{2}+\|(-\Delta)^{\frac{\alpha}{2}}\theta(t)\|^{2}_{L^{2}(\R^{d})}=0.
\end{equation}
Note that this estimate does not involve any norms of the velocity at all.  However, to prove a H\"older continuity result, one needs to understand the balance of energy  within a small ball.  In Section \ref{S.Cacciop}, we study the evolution of the quantity $t \to \|\eta^{2}\theta(t)\|_{L^{2}(\R^{d})}$ for an appropriate cutoff $\eta$.  The contribution of the advection $u \cdot \nabla \theta$ requires one to estimate 
\begin{equation}\label{eq:TransportContr}
2\left|\left\langle u \cdot \nabla \theta, \eta^{2}\theta \right\rangle\right| = \left|\left\langle u, \theta^2 \nabla(\eta^{2}) \right\rangle\right|.
\end{equation}
The main result of this section is Lemma \ref{L.ConditionalNonlinearIneq}, which gives an estimate for \eqref{eq:TransportContr} in terms of $u$ in $\text{BMO}^{-\gamma}(\R^{d})$ and localized norms of $\theta$ and $\theta^{*}$.

\subsection{Classical auxilliary results}\label{S.AuxillaryResults}
In this subsection, we collect several auxilliary results about the fractional laplace operator and the space $\text{BMO}(\R^{d})$ that will be used in the proofs of Lemmas \ref{L.BMOProd} and \ref{L.DriftEstimate} below.  The first result, the Kato-Ponce inequality, compensates for the lack of a pointwise product rule for the fractional Laplacian.
\begin{Thm}\label{T.KatoPonce}
Given $\gamma>0$, $p\in [1,\infty)$ and exponents $q,r \in (1,\infty)$ such that $\frac{1}{r}+\frac{1}{q}=\frac{1}{p}$, there exists a constant $C=C(d,\gamma,q,r,p)$ such that for all $f,g\in W^{\gamma,q}(\R^d)\cap L^r(\R^d)$,
\begin{align}\label{eq:KatoPonce}
&\big \|(-\Delta)^{\frac{\gamma}{2}}(fg) \big \|_{L^{p}(\R^{d})}
\leq C \left[ \big\|(-\Delta)^{\frac{\gamma}{2}}f \big\|_{L^{q}(\R^{d})}\|g\|_{L^{r}(\R^{d})} 
+\big\|(-\Delta)^{\frac{\gamma}{2}}g\big\|_{L^{q}(\R^{d})}\|f\|_{L^{r}(\R^{d})} \right].
\end{align}
\end{Thm}
\begin{proof}
Theorem 1 in \cite{grafakos2014kato} proves the same inequality for $f,\;g\in \mathcal{S}(\R^d)$. The version of the Theorem stated above can be proved using density of $\mathcal{S}(\R^d)$ in $W^{\gamma,q}(\R^d)\cap L^r(\R^d)$.  Clearly one can pass to the limit on the RHS of \eqref{eq:KatoPonce}. To treat the LHS, note that if 
$f_n\to f$ and $g_n\to g$ in \linebreak $W^{\gamma,q}(\R^d)\cap L^r(\R^d)$, then $(-\Delta)^{\frac{\gamma}{2}}(f_ng_n)\to (-\Delta)^{\frac{\gamma}{2}}(fg)$ in $\mathcal{S}'(\R^d)$.  By the lower-semicontinuity of the $L^p(\R^d)$ norm, we obtain the desired inequality in the limit.
\end{proof}
The basic mechanism for using the regularity of $\theta$ to compensate for the irregularity of $u$ is the the fractional integration by parts formula below.
\begin{Lem} \label{L.FractionalIBP}
For all $f \in C^{\infty}(\R^{d}) \cap \text{BMO}(\R^{d})$ and $g \in W_c^{1,\infty}(\R^{d})$,
\begin{equation}\label{eq:FracIBP}
\int_{\R^{d}}(-\Delta)^{\frac{\gamma}{2}}f g\; \dx = \int_{\R^{d}} f (-\Delta)^{\frac{\gamma}{2}}g\; \dx.
\end{equation}
\end{Lem}
Observe that if $g$ in Lemma \ref{L.FractionalIBP} is compactly supported but $f$ is not, then fractional integration by parts turns a localized integral into an integration over the whole space.  The next Lemma allows us to use the decay of $(-\Delta)^{\frac{\gamma}{2}}g$ away from the support of $g$ to estimate the far field contribution to this integral.  Also, it allows us to freely subtract a constant from $f$, which will be useful when applying the John-Nirenberg inequality.
\begin{Lem} \label{L.AwayFromSupport}
For each $\gamma \in (0,1)$, there exists a positive constant $C_{\gamma}$ with the following property.  For all $g \in W_c^{1,\infty}(\R^{d})$ supported in $B_{1}$ and each $x \in B_{2}^{c}$,
\begin{align}
\big|(-\Delta)^{\frac{\gamma}{2}}g(x)\big| \leq C_{\gamma}\|g\|_{L^{1}(B_{1})}\big ( 1+|x|^{d+\gamma} \big )^{-1}.
\end{align}
In particular, $(-\Delta)^{\frac{\gamma}{2}}g \in L^{1}(\R^{d})$ and satisfies:
\begin{equation}\label{eq:MeanZero}
\int_{\R^{d}}(-\Delta)^{\frac{\gamma}{2}}g(x)\dx = 0.
\end{equation}
\end{Lem}
\begin{proof}
Let $x \in B_{2}^{c}$, then since $g$ is compactly supported in $B_{1}$
\begin{equation} \label{eq:fracLapAwayFromSupp}
(-\Delta)^{\frac{\gamma}{2}}g(x)=c_{\gamma}\int_{B_{1}}\frac{g(y)}{|x-y|^{d+\gamma}}\dee y,
\end{equation}
where the integral is absolutely convergent.  Moreover, for $y \in B_{1}$ we have the elementary inequality
\begin{equation} \label{eq:ElemIneq}
|x-y| \geq |x|-|y|=\frac{1}{3}|x|+\bigg(\frac{2}{3}|x|-|y|\bigg) \geq \frac{1}{3}(1+|x|).
\end{equation}
Combining \eqref{eq:fracLapAwayFromSupp} and \eqref{eq:ElemIneq} gives the claim.

To show \eqref{eq:MeanZero}, we use again the integral definition of the fractional Laplacian. As $g\in W^{1,\infty}_c(\R^d)$ and $\gamma\in (0,1)$, we may omit the principal value part of the definition. Therefore,
$$
\int_{\R^d} (-\Delta)^{\frac{\gamma}{2}}g(x)\;\dx=\int_{\R^d} \int_{\R^d}\frac{g(x)-g(y)}{|x-y|^{d+\gamma}}\;\dx\dee y=0,
$$
due to the anti-symmetry of the integrand.
\end{proof}
Finally, we need some classical estimates for $\text{BMO}(\R^{d})$ functions, a local control through the John-Nirenberg inequality, as well as an estimate on the far field behavior.  
\begin{Prop} \label{P.ClassicalBMO} 
For all $p \in (1,\infty)$, there exists a positive constant $C_{p}$ such that for each \linebreak $\psi \in \text{BMO}(\R^{d})$,
\begin{equation}\label{eq:John_Nirenberg}
\|\psi-M(\psi,B_{2})\|_{L^{p}(B_{2})}\leq C_{p}\|\psi\|_{\text{BMO}(B_{2})}.
\end{equation} 
Moreover, for all $\gamma \in (0,2)$, there exists a positive constant $C_{\gamma}$ such that:
\begin{equation}\label{eq:FarFieldGrowth}
\int_{\R^{d}}\frac{|\psi-M(\psi,B_{2})|}{2^{d+\gamma}+|x|^{d+\gamma}}\dx \leq C_{\gamma}\|\psi\|_{\text{BMO}(\R^{d})},
\end{equation}
where
$$
M(\psi,B_{2})=\frac{1}{|B_{2}|}\int_{B_2}\psi \dx.
$$
\end{Prop}
A proof of \eqref{eq:John_Nirenberg} may be found in \cite{stein2016harmonic} and the proof of \eqref{eq:FarFieldGrowth} appears in \cite{strichartz1980bounded}.


\subsection{A global regularity trade-off}
Now we employ the auxilliary results of Section \ref{S.AuxillaryResults} to prove the following:
\begin{Lem}\label{L.BMOProd} 
For all $p \in (1,\infty)$ and $\gamma \in (0,1)$, there exists a positive constant $C_{p,\gamma}$ such that the following is true.  For each $f \in C^{\infty}(\R^{d}) \cap \text{BMO}^{-\gamma}(\R^{d})$ and $g \in W_c^{1,\infty}(\R^{d})$ supported in $B_{1}$,
\begin{equation}\label{E.BMOProd}
\int_{\R^{d}}fg\dx \leq C_{p,\gamma}\|f\|_{ \text{BMO}^{-\gamma}(\R^{d})}\bigg [ \big\|(-\Delta)^{\frac{\gamma}{2}}g \big\|_{L^{p}(\R^{d})}+\|g\|_{L^{1}(\R^{d})}\bigg].
\end{equation}
\end{Lem}
\begin{proof}\,Begin by writing $f=(-\Delta)^{\frac{\gamma}{2}}\psi$ and applying Lemma \ref{L.FractionalIBP} to find:
$$
\int_{\R^{d}}fg\dx = \int_{\R^{d}}\psi (-\Delta)^{\frac{\gamma}{2}}g\dx.
$$
By Lemma \ref{L.AwayFromSupport}, $(-\Delta)^{ \frac{\gamma}{2} }g$ integrates to zero over $\R^{d}$.  We may exploit this by introducing:
$$
\psi_{0}(x)=\psi(x)-M(\psi,B_{2}).
$$
A free subtraction yields:
\begin{align}
\int_{\R^{d}}fg\dx = \int_{\R^{d}}\psi_{0} (-\Delta)^{\frac{\gamma}{2}}g\dx 
=\int_{B_{2}}\psi_{0} (-\Delta)^{\frac{\gamma}{2}}g\dx+\int_{B_{2}^{c}}\psi_{0} (-\Delta)^{\frac{\gamma}{2}}g\dx.
\end{align}
To estimate the inner contribution, apply H\"{o}lder's inequality followed by Proposition \ref{P.ClassicalBMO} to obtain: 
\begin{align}
\int_{B_{2}}\psi_{0} (-\Delta)^{\frac{\gamma}{2}}g\dx &\leq \|\psi_{0}\|_{L^{\frac{p}{p-1}}(B_{2})}\big \|(-\Delta)^{\frac{\gamma}{2}}g\big \|_{L^{p}(\R^{d})}\\
& \leq C_{p}\|\psi\|_{\text{BMO}(\R^{d})}\big \|(-\Delta)^{\frac{\gamma}{2}}g\big \|_{L^{p}(\R^{d})}. 
\end{align}
To estimate the outer contribution, use Lemma \ref{L.AwayFromSupport} (noting that $g$ is compactly supported in $B_{1}$) followed by Proposition \ref{P.ClassicalBMO} to deduce
\begin{align}
\int_{B_{2}^{c}}\psi_{0} (-\Delta)^{\frac{\gamma}{2}}g\dx &\leq C_{\gamma}\|g\|_{L^{1}(B_{1})}  \int_{\R^{d}}\frac{|\psi(x)-M(\psi,B_{2})|}{1+|x|^{d+\gamma}}\dx \\
&\leq C_{\gamma} \|g\|_{L^{1}(B_{1})} \int_{\R^{d}}\frac{|\psi(x)-M(\psi,B_{2})|}{ 2^{d+\gamma}+|x|^{d+\gamma}}\dx \\
&\leq C_{\gamma}\|g\|_{L^{1}(B_{1})}\|\psi\|_{\text{BMO}(\R^{d})}.
\end{align}
The proof of the Lemma is completed in view of the observation $\|\psi\|_{\text{BMO}(\R^{d})}=\|f\|_{\text{BMO}^{-\gamma}(\R^{d})}$.
\end{proof}


\subsection{A local regularity trade-off}
We are now in a position to prove the main result of this section.
\begin{Lem}\label{L.DriftEstimate}
For each $\gamma \in (0,\alpha)$, there exists a positive constant $C_{\gamma,\alpha}$ with the following property.  For all $u \in C^{\infty}(\R^{d}) \cap BMO^{-\gamma}(\R^{d})$, $\theta \in C^{\infty}(\R^{d})$, and $\sigma \in C^{\infty}_{c}\big(B_1\times [0,\infty) \big)$, the inequality below holds true:
\begin{align}
\left | \int_{\R^{d}}u \cdot \eta \nabla \eta \theta_{+}^{2}\dee x \right |
&\leq C_{\gamma,\alpha}\nu^{-1}\left[1+\|u\|_{\text{BMO}^{-\gamma}(\R^{d})}^{2}\right]\|(\eta + |\nabla \eta|)\theta_{+}\|^{2}_{L^{2}(\R^{d})}\\
&+\|\nabla(\sigma \theta^{*}_{+})\|_{L^{2}_{\omega_{\alpha}}(\R^{d} \times \R_{+})}+\nu \|\nabla(\nabla_{x} \sigma \theta_{+}^{*})\|^{2}_{L^{2}_{\omega_{\alpha}}(\R^{d} \times \R_{+})},
\end{align}
where $\nu \in (0,1)$ may be arbitrary, $\eta(x)=\sigma(x,0),$ and $\theta_{+}^{*}=\left(\theta^{*} \right)_{+}$.
\end{Lem}

\begin{proof}
Let $f=u$ and $g=\eta \nabla \eta\theta_{+}^{2}$.  Since $\theta \in C^{\infty}(\R^{d})$, it follows that $g \in W_c^{1,\infty}(\R^{d})$.  Hence, we may apply Lemma \ref{L.BMOProd} with an exponent $p(\alpha,\gamma) \in (1,2)$ to be chosen below.  This yields:
\begin{equation}
\int_{\R^{d}}u \cdot \eta \nabla \eta \theta_{+}^{2}\dx \leq C_{\gamma,\alpha}\|u\|_{\text{BMO}^{-\gamma}(\R^{d})}\bigg [ \|\eta \nabla \eta \theta_{+}^{2}\|_{L^{1}(\R^{d})}+\big \|(-\Delta)^{\frac{\gamma}{2}}(\eta \nabla \eta \theta_{+}^{2})\big \|_{L^{p}(\R^{d})} \bigg].
\end{equation}
To estimate the first term, use Young's inequality to find:
\begin{equation}
\|u\|_{BMO^{-\gamma}(\R^{d})} \|\eta \nabla \eta \theta_{+}^{2}\|_{L^{1}(\R^{d})}
\leq \frac{1}{2}\left[1+\|u\|_{\text{BMO}^{-\gamma}(\R^{d})}^{2} \right ]\|(\eta + |\nabla \eta|)\theta_{+}\|_{L^{2}(\R^{d})}^{2}.
\end{equation}
The next step is to apply Kato-Ponce's inequality (Theorem \ref{T.KatoPonce}) to $\eta \theta_{+}$ and $\nabla \eta \theta_{+}$ to obtain:
\begin{align}
&\big\|(-\Delta)^{\frac{\gamma}{2}}(\eta \theta_{+} \nabla \eta \theta_{+})\big\|_{L^{p}(\R^{d})}\\
&\leq C_{\gamma}\bigg [ \big\|(-\Delta)^{\frac{\gamma}{2}}(\eta \theta_{+})\big\|_{L^{\frac{2p}{2-p}}(\R^{d})}\|\nabla \eta \theta_{+}\|_{L^{2}(\R^{d})}+\|\eta \theta_{+}\|_{L^{2}(\R^{d})}\big\|(-\Delta)^{\frac{\gamma}{2}}(\nabla \eta \theta_{+})\big \|_{L^{\frac{2p}{2-p}}(\R^{d})} \bigg].
\end{align}
Since $\alpha>\gamma$, we may choose the exponent $p=d/(d-\alpha+\gamma)$ which is less than $2$ since $\gamma>0$, $\alpha<1$ and $d \geq 2$. This allows us to use the Hardy-Littlewood-Sobolev inequality to deduce 
\begin{align}
&\big\|(-\Delta)^{\frac{\gamma}{2}}(\eta \theta_{+} \nabla \eta \theta_{+})\big\|_{L^{p}(\R^{d})}\\
&\leq C_{\gamma}\bigg [ \big\|(-\Delta)^{\frac{\alpha}{2}}(\eta \theta_{+})\big \|_{L^{2}(\R^{d})}\|\nabla \eta \theta_{+}\|_{L^{2}(\R^{d})}+\|\eta \theta_{+}\|_{L^{2}(\R^{d})}\big \|(-\Delta)^{\frac{\alpha}{2}}(\nabla \eta \theta_{+})\big \|_{L^{2}(\R^{d})} \bigg].
\end{align}
Combining these observations and using Young's inequality with a parameter $\nu>0$, we find:
\begin{align}
\left | \int_{\R^{d}}u \cdot \eta \nabla \eta \theta_{+}^{2}\dx \right | &\leq C_{\gamma}\nu^{-1}\bigg[1+\|u\|_{\text{BMO}^{-\gamma}(\R^{d})}^{2} \bigg ]\|(\eta + |\nabla \eta|)\theta_{+}\|^{2}_{L^{2}(\R^{d})}\\
&+\big \|(-\Delta)^{\frac{\alpha}{2}}(\eta \theta_{+})\big \|^{2}_{L^{2}(\R^{d})}+\nu \big \|(-\Delta)^{\frac{\alpha}{2}}(\nabla \eta \theta_{+})\big \|^{2}_{L^{2}(\R^{d})}.
\end{align}
The proof may be completed by using the variational characterization of the $\dot{H}^{\alpha}$ norm, Proposition \ref{prop:variational}.  Namely, choose $\sigma \theta_{+}^{*}$ to extend $\eta \theta_{+}$ and $\nabla_{x} \sigma \theta_{+}^{*}$ to extend $\nabla \eta \theta_{+}$.
\end{proof}

\section{A Caccioppoli Type Inequality}\label{S.Cacciop}
In this section, we use Lemma \ref{L.DriftEstimate} to obtain a form of Caccioppoli's inequality. A key ingredient is Lemma \ref{L.FundamentalIdentity}, which allows us to lift a fractional derivative operator (for which the product rule does not hold) to a local derivative operator of the Caffarelli-Silvestre extension (by adding a dimension).  On the extension, the usual Caccioppoli manipulations for local elliptic/parabolic equations can be performed.  The result is an inequality involving both $\theta$ and $\theta^{*}$, which will have a cost when we perform the De Giorgi iteration in Sections \ref{S.OscRed} and \ref{S.Isoperim}.

The main result of this section is Proposition \ref{C.NoCutoffCacciopioli}, whose proof relies on the following Lemma. 

\begin{Lem} \label{P.CacciopioliInequality}
There exists a positive constant $C(D,\alpha)$ such that for all $\theta \in \mathcal{S}(D,\alpha)$, the following property is true.  For each $\sigma \in C^{\infty}_{c}\big (\R^{d}\times [0,\infty) \big )$ and $0<s<t\le 6$ $:$
\begin{equation}\label{eq:CacciopStep1}
\begin{split}
\|\eta \theta_{+}(t)\|_{L^{2}(\R^{d})}^{2}&+\int_{s}^{t}\|\nabla \big (\sigma \theta_{+}^{*}(r) \big )\|^{2}_{L_{\omega_\alpha}^{2}(\R^d\times\R_+)}\dee r \leq \|\eta \theta_{+}(s)\|_{L^{2}(\R^{d})}^{2}\\
&+ 3\int_{s}^{t}\|\nabla \sigma \theta_{+}^{*}(r)\|^{2}_{L_{\omega_\alpha}^{2}(\R^d\times\R_+)}\dee r
+C\nu^{-1}\left[ \int_{s}^{t}\|(\eta+|\nabla \eta|)\theta_{+}(r)\|_{L^{2}(\R^{d})}^{\frac{2\alpha}{1-\alpha}}\dee r \right]^{\frac{1}{\alpha}-1}\\
&+\nu \int_{s}^{t}\|\nabla \big (\nabla_{x} \sigma \theta_{+}^{*}(r) \big )\|^{2}_{L_{\omega_\alpha}^{2}(\R^d\times\R_+)}\dee r, \\
\end{split}
\end{equation}
where the constant $\nu \in (0,1)$ may be chosen arbitrarily and $\eta(x)=\sigma(x,0)$.
\end{Lem}

\begin{proof}
Multiplying the equation for $\theta \in \mathcal{S}(D,\alpha)$ by $\eta^{2}\theta_{+}$ and integrating over $[s,t] \times \R^{d}$ yields:
\begin{equation} \label{E.LocalizeIBP}
\begin{split}
&\frac{1}{2}\int_{\R^{d}}\theta_{+}^{2}(t,x)\eta^{2}(x)\dx + \int_{s}^{t}\int_{\R^{d}} (-\Delta)^{\alpha}\theta \eta^{2}\theta_{+} \dx\dee r\\
&= \frac{1}{2}\int_{\R^{d}}\theta_{+}^{2}(s,x)\eta^{2}(x)\dx - \frac{1}{2}\int_{s}^{t}\int_{\R^{d}} u \cdot \nabla(\theta_{+}^{2}) \eta^{2} \dx\dee r.
\end{split}
\end{equation}
Since $\Div u = 0$, we may write:
\begin{align}\label{E.DivFree}
&- \frac{1}{2}\int_{s}^{t}\int_{\R^{d}} u \cdot \nabla(\theta_{+}^{2}) \eta^{2} \dx\dee r
=\frac{1}{2}\int_{s}^{t}\int_{\R^{d}} u \cdot \nabla(\eta^{2})\theta_{+}^{2}\dx\dee r.
\end{align}
Next we will verify that for all times in $[s,t]$,
\begin{equation}\label{E.LiftToUHP}
\int_{\R^{d}} (-\Delta)^{\alpha}\theta \eta^{2} \theta_{+} \dee x = \int_{\R^{d}}\int_{0}^{\infty}z^{1-2\alpha}\bigg ( |\nabla(\sigma\theta_{+}^{*})|^{2}-|\nabla \sigma|^{2}|\theta_{+}^{*}|^{2} \bigg)\dx \dee z.
\end{equation}
Applying Lemma \ref{L.FundamentalIdentity} with $f = \theta$, $g_{0} = \eta^{2} \theta_{+}$, and the extension $g=\sigma^{2}\theta_{+}^{*}$,
\begin{align}
\int_{\R^{d}} (-\Delta)^{\alpha}\theta \eta^{2} \theta_{+} \dx &= \int_{0}^{\infty}\int_{\R^{d}} z^{1-2\alpha}\nabla(\sigma^{2}\theta_{+}^{*}) \cdot \nabla \theta^{*}\dx\dee z \\
& = \int_{0}^{\infty}\int_{\R^{d}} z^{1-2\alpha}\nabla(\sigma^{2}\theta_{+}^{*}) \cdot \nabla \theta_{+}^{*}\dx\dee z.
\end{align}
The identity \eqref{E.LiftToUHP} now follows from the relation 
$$
\nabla(\sigma^{2}\theta_{+}^{*}) \cdot \nabla \theta^{*}_{+} = |\nabla(\sigma \theta^{*}_{+})|^{2}-|\nabla \sigma|^{2}|\theta^{*}_{+}|^{2}.
$$
Substituting the identities \eqref{E.DivFree} and \eqref{E.LiftToUHP} into \eqref{E.LocalizeIBP} and multiplying by a factor of $2$ yields:
\begin{equation}\label{E.ReductionToDrift}
\begin{split}
\|\eta \theta_{+}(t)\|_{L^{2}(\R^{d})}^{2}&+2\int_{s}^{t}\|\nabla \big (\sigma \theta_{+}^{*}(r) \big )\|^{2}_{L^{2}_{\omega_{\alpha}}(\R^{d} \times \R_{+})}\dee r
\leq \|\eta \theta_{+}(s)\|_{L^{2}(\R^{d})}^{2} \\
&+ 2\int_{s}^{t}\|\nabla \sigma \theta_{+}^{*}(r)\|^{2}_{L^{2}_{\omega_{\alpha}}(\R^{d} \times \R_{+})}\dee r+\frac{1}{2}\int_{s}^{t}\int_{\R^{d}} u \cdot \nabla(\eta^{2})\theta_{+}^{2}\dx\dee r.
\end{split}
\end{equation}
Fix a $\nu>0$ and apply Lemma \ref{L.DriftEstimate} for each time $r \in (s,t)$, then integrate to find
\begin{equation} \label{E.DriftEstCacc}
\begin{aligned}
\int_{s}^{t}\int_{\R^{d}}u \cdot \eta \nabla \eta \theta_{+}^{2}\dx \dee r &\leq C_{\gamma}\nu^{-1}\int_{s}^{t}\left[1+\|u(r)\|_{\text{BMO}^{-\gamma}(\R^{d})}^{2} \right ]\|(\eta + |\nabla \eta|)\theta_{+}(r)\|^{2}_{L^{2}(\R^{d})}\dee r\\
&+\int_{s}^{t}\|\nabla(\sigma \theta^{*}_{+})\|^{2}_{L^{2}_{\omega_{\alpha}}(\R^{d} \times \R_{+})}\dee r+\nu \int_{s}^{t}\|\nabla(\nabla_{x} \sigma \theta_{+}^{*})\|^{2}_{L^{2}_{\omega_{\alpha}}(\R^{d} \times \R_{+})}\dee r.
\end{aligned}
\end{equation}
Note that the second integral on the RHS of \eqref{E.DriftEstCacc} can be absorbed into the LHS of \eqref{E.ReductionToDrift}.  Hence, to complete the proof it suffices to estimate the first term on the RHS of \eqref{E.DriftEstCacc}. Applying Holder's inequality with the exponent $\alpha/(2\alpha-1)$ yields:
\begin{align}
&\int_{s}^{t}\left[1+\|u(r)\|_{\text{BMO}^{-\gamma}(\R^{d})}^{2}\right]\|(\eta + |\nabla \eta|)\theta_{+}(r)\|^{2}_{L^{2}(\R^{d})}\dee r \\ 
&\leq \bigg ( \int_{s}^{t}\left[1+\|u(r)\|_{\text{BMO}^{-\gamma}(\R^{d})}^{2}\right]^{\frac{\alpha}{2\alpha-1}}\dee r \bigg )^{2-\frac{1}{\alpha}} \bigg ( \int_{s}^{t}\|(\eta + |\nabla \eta|)\theta_{+}(r)\|^{\frac{2\alpha}{1-\alpha}}_{L^{2}(\R^{d})}\dee r \bigg )^{\frac{1}{\alpha}-1}
\end{align}
Recall that $\theta \in \mathcal{S}(D,\alpha)$ implies $u \in \LB{q}{-\gamma}$, where the exponents $q$ and $\gamma$ satisfy $\frac{2\alpha}{q}+\gamma=2\alpha-1$.  Since $\gamma>0$, it follows that $q>2\alpha/(2\alpha-1)$, so applying Holder's inequality again,
\begin{align}
\bigg ( \int_{s}^{t}\left[1+\|u(r)\|_{\text{BMO}^{-\gamma}(\R^{d})}^{2}\right]^{\frac{\alpha}{2\alpha-1}}\dee r \bigg )^{2-\frac{1}{\alpha}} 
&\leq C_{\alpha}\left [1+\|u\|_{\LB{\frac{2\alpha}{2\alpha-1}}{-\gamma}}^{2}\right]\\
&\leq C_{\alpha}\left [1+\|u\|_{\LB{q}{-\gamma}}^{2} \right] \\
&\leq C(D,\alpha).
\end{align}
Note that we used the fact that $s,t\leq 6$ in the second inequality.  Combining these observations with \eqref{E.ReductionToDrift} and \eqref{E.DriftEstCacc} completes the proof.

\end{proof}

In order for \eqref{eq:CacciopStep1} to be useful in the De Giorgi iteration step, we need to remove the last term on the RHS.  The strategy is to do this by iterating the estimate.    The first step is to choose particular cutoffs, yielding an inequality similar to \eqref{eq:CacciopStep1}, but more amenable to the hypotheses of Lemma \ref{L.IterationLemma}.  The end result is that we remove the overlap at the expense of the precision of the estimate.  The details are given in the proof of the next Proposition. 
\begin{Prop} \label{C.NoCutoffCacciopioli}
There exists a positive universal constant $C(D,\alpha)$ such that for all $\theta \in \mathcal{S}(D,\alpha)$, the following inequality holds.  For all cube sizes $0<r_{1}<r_{2}\leq 3$, heights $0<\delta_{1}<\delta_{2}\leq 3$, and times $0<t_{1}<t_{2}\leq 6$ $:$

\begin{equation} \label{eq:NoCutoffCacciopoli}
\begin{split}
\|\theta_{+}(t_{2})\|_{L^{2}(B_{r_{1}})}^{2} &+ C^{-1}\int_{t_1}^{t_2}\|\nabla \theta_{+}^{*}\|^{2}_{L^{2}_{\omega_\alpha}(B_{r_{1}} \times (0,\delta_{1}))}\;\dt 
\leq
\|\theta_{+}(t_{1})\|_{L^{2}(B_{r_{2}})}^{2}\\
&+\frac{C}{(r_{2}-r_{1})^{4}(\delta_{2}-\delta_{1})^{2}}\left [ \int_{t_1}^{t_2}\|\theta_{+}^{*}\|^{2}_{L^{2}_{\omega_\alpha}(B_{r_{2}} \times [0,\delta_{2}])}\;\dt +\left(\int_{t_1}^{t_2}\|\theta_{+}\|^{\frac{\alpha}{1-\alpha}}_{L^{2}(B_{r_{2}})}\;\dt\right)^{\frac{1-\alpha}{\alpha}}\right ]. 
\end{split}
\end{equation} 
\end{Prop}
\begin{proof} Fix $r,\;R,\;\delta$ and $\rho$, such that $0<r_1<r<R<r_2 \leq 3$ and $0<\delta_1<\delta<\rho<\delta_2 \leq 3$. Also, fix a cutoff function $\sigma \in C_{c}^{\infty}(\R^{d} \times [0,\infty))$ such that:
\begin{equation}
1_{B_{r} \times [0,\delta]} \leq \sigma \leq 1_{B_{R} \times [0,\rho]},
\end{equation}
and the derivatives of the cutoff satisfy:
\begin{align}
\|\nabla \sigma\|_{L^{\infty}(\R^{d} \times \R_{+})}
\leq C(R-r)^{-1}(\rho-\delta)^{-1}, \quad
\|\nabla^{2}\sigma\|_{L^{\infty}(\R^{d} \times \R_{+})}&\leq C(R-r)^{-2}(\rho-\delta)^{-2}.
\end{align}
Using Proposition \ref{P.CacciopioliInequality} and noting the support of $\sigma$ yields for each $\nu \in (0,1)$:
\begin{equation}\label{E.ApplyCacciop}
\begin{split}
\|\theta_{+}(t_{2})\|_{L^{2}(B_{r})}^{2}&+\int_{t_{1}}^{t_{2}}\big\|\nabla \big(\theta_{+}^{*}(t)\big)\big\|^{2}_{L_{\omega_{\alpha}}^{2}(B_{r}\times [0,\delta])}\dt
\leq \|\theta_{+}(t_{1})\|_{L^{2}(B_{R})}^{2} \\
&+ 3\int_{t_{1}}^{t_{2}}\|\nabla \sigma \theta_{+}^{*}(t)\|^{2}_{L_{\omega_{\alpha}}^{2}(\R^{d} \times \R_{+})}\dt
+\nu \int_{t_{1}}^{t_{2}}\big\|\nabla \big(\nabla \sigma \theta_{+}^{*}(t) \big ) \big\|^{2}_{L_{\omega_{\alpha}}^{2}(\R^{d} \times \R_{+})}\dt\\
&+C\nu^{-1}\left( \int_{t_{1}}^{t_{2}}\|(\eta+|\nabla \eta|)\theta_{+}(t)\|_{L^{2}(\R^{d})}^{\frac{2\alpha}{1-\alpha}}\dt \right)^{\frac{1}{\alpha}-1}.
\end{split}
\end{equation}
Observe that for all times in $t \in [t_{1},t_{2}]$:
\begin{align}
&\big\|\nabla \big (\nabla \sigma \theta_{+}^{*}(t) \big )\big\|^{2}_{L_{\omega_{\alpha}}^{2}(\R^{d}\times \R_{+})}\\
&\leq 2\|\nabla^{2}\sigma\|^{2}_{L^{\infty}(\R^{d} \times \R_{+})}\|\theta_{+}^{*}(t)\|^{2}_{L_{\omega_{\alpha}}^{2}(B_{R} \times (0,\rho))}
+
2\|\nabla \sigma\|_{L^{\infty}(\R^{d} \times \R_{+})}^{2}\big\|\nabla \theta_{+}^{*}(t)\big\|^{2}_{L_{\omega_{\alpha}}^{2}(B_{R} \times (0,\rho))}\\
&\leq 2C(\rho-\delta)^{-4}(R-r)^{-4}\|\theta_{+}^{*}(t)\|^{2}_{L_{\omega_{\alpha}}^{2}(B_{R} \times (0,\rho))}+2C(\rho-\delta)^{-2}(R-r)^{-2}\big\|\nabla \theta_{+}^{*}(t)\big\|^{2}_{L_{\omega_{\alpha}}^{2}(B_{R} \times (0,\rho))}.
\end{align}
Hence, choosing $\nu = (2C)^{-1}(\rho-\delta)^{2}(R-r)^{2}$ yields:
\begin{equation}\label{E.OverlapTerms}
\begin{split}
& 3\int_{t_1}^{t_2}\|\nabla \sigma \theta_{+}^{*}\|^{2}_{L_{\omega_{\alpha}}^{2}(\R^{d} \times \R_{+})}\dee t
+\nu \int_{t_1}^{t_2}\big\|\nabla \big (\nabla \sigma \theta_{+}^{*} \big )\big\|^{2}_{L_{\omega_\alpha}^{2}(\R^{d} \times \R_{+})}\dee t\\
& \leq \frac{1}{2}\int_{t_1}^{t_2}\|\nabla \theta_{+}^{*}\|^{2}_{L_{\omega_\alpha}^{2}(B_{R}\times (0,\rho))}\;\dee t 
+
C(\rho-\delta)^{-2}(R-r)^{-2}\int_{t_1}^{t_2}\|\theta_{+}^{*}\|^{2}_{L_{\omega_\alpha}^{2}(B_{R}\times(0,\rho))}\;\dee t\\
& \leq \frac{1}{2}\int_{t_1}^{t_2}\|\nabla \theta_{+}^{*}\|^{2}_{L_{\omega_\alpha}^{2}(B_{R}\times (0,\rho))}\;\dee t
+
C(\rho-\delta)^{-2}(R-r)^{-4}\int_{t_1}^{t_2}\|\theta_{+}^{*}\|^{2}_{L_{\omega_\alpha}^{2}(B_{R}\times(0,\rho))}\;\dee t.
\end{split}
\end{equation}
Moreover, using the bounds on the gradient of the cutoff
\begin{equation}\label{E.UsualTerms}
C\nu^{-1}\left ( \int_{t_1}^{t_2}\|(\eta+|\nabla \eta|)\theta_{+}\|_{L^{2}(\R^{d})}^{\frac{2\alpha}{1-\alpha}}\dee t \right )^{\frac{1-\alpha}{\alpha}}
 \leq C(\rho-\delta)^{-2}(R-r)^{-4}\left(\int_{t_1}^{t_2}\|\theta_{+}\|^{\frac{\alpha}{1-\alpha}}_{L^{2}(B_{R})}\;\dee t\right)^{\frac{1-\alpha}{\alpha}}.
\end{equation}
Using the inequalities \eqref{E.OverlapTerms} and \eqref{E.UsualTerms} to update the estimate \eqref{E.ApplyCacciop} yields:
\begin{equation} 
\begin{split}
\|\theta_{+}(t_{2})\|_{L^{2}(B_{r})}^{2} &+ \int_{t_1}^{t_2}\|\nabla \theta_{+}^{*}\|^{2}_{L^{2}_{\omega_\alpha}(B_{r} \times (0,\delta))}\;\dt \leq \|\theta_{+}(t_{1})\|_{L^{2}(B_{R})}^{2}\\
&+\frac{C}{(R-r)^{4}(\rho-\delta)^{2}}\left [ \int_{t_1}^{t_2}\|\theta_{+}^{*}\|^{2}_{L^{2}_{\omega_\alpha}(B_{R} \times [0,\rho])}\;\dt +\left(\int_{t_1}^{t_2}\|\theta_{+}\|^{\frac{\alpha}{1-\alpha}}_{L^{2}(B_{R})}\;\dt\right)^{\frac{1-\alpha}{\alpha}}\right ]\\
&+\frac{1}{2}\int_{t_1}^{t_2}\|\nabla \theta_{+}^{*}\|^{2}_{L_{\omega_\alpha}^{2}(B_{R}\times (0,\rho))}\;\dee t. 
\end{split}
\end{equation} 

Finally, we just need to apply this inequality iteratively to obtain the desired result. This iteration procedure is encoded in Lemma~\ref{L.IterationLemma}, where we use with the following parameters:
\begin{align}
&A_{1}=\big\|\theta_{+}(t_{1})\big\|_{L^{2}(B_{r_{2}})}^{2}-\big\|\theta_{+}(t_{2})\big\|_{L^{2}(B_{r_{1}})}^{2}\\
&A_{2}=C \left [  \int_{t_1}^{t_2}\|\theta_{+}^{*}\|^{2}_{L^{2}_{\omega_\alpha}(B_{r_{2}} \times [0,\delta_{2}])} +\left(\int_{t_1}^{t_2}\|\theta_{+}\|^{\frac{\alpha}{1-\alpha}}_{L^{2}(B_{r_{2}})}\;\dt\right)^{\frac{1-\alpha}{\alpha}}
\right ]\\
& \phi(r,\delta)=\int_{t_1}^{t_2}\|\nabla \theta_{+}^{*}\|^{2}_{L_{\omega_\alpha}^{2}(B_{r}\times (0,\delta))}\;\dee t \\
& \kappa =\frac{1}{2}, \; \beta_{1}=2, \; \beta_{2}=4.
\end{align}
As the constant obtained from iterating depends only on $\kappa$, we deduce the desired inequality.
\end{proof}

\begin{Lem}\label{L.IterationLemma}
Let $\phi: [r_{1},r_{2}] \times [\delta_{1},\delta_{2}] \to \R_{+}$ be a positive, bounded function which is monotone in the two parameter sense.  Let $A_{1} \in \R$, $A_{2} \in \R_{+}$ be constants and $\beta_{1},\beta_{2} \in \R_{+}$ be exponents. 

Assume there exists a $\kappa \in (0,1)$ with the following property.  For all $r_{1}<r<R<r_{2}$ and $\delta_{1}<\delta<\rho<\delta_{2}$ $:$
\begin{equation}
\phi(r,\delta) \leq A_{1}+\frac{A_{2}}{(\delta-\rho)^{\beta_{1}}(R-r)^{\beta_{2}}}+\kappa \phi(R,\rho).
\end{equation}
Then there exists a constant $C(\kappa)$ depending solely on $\kappa$, such that:
\begin{equation}
\phi(r_{1},\delta_{1}) \leq C(\kappa) \left( A_{1}+\frac{A_{2}}{(\delta_{2}-\delta_{1})^{\beta_{1}}(r_{2}-r_{1})^{\beta_{2}}} \right).
\end{equation}
\end{Lem}

\section{Oscillation Reduction}\label{S.OscRed}
In the classical De Giorgi theory of elliptic equations in divergence form, the first step towards H\"older continuity is to go from $L^{2}$ to $L^{\infty}$.  Namely, one proves that the $L^{2}$ norm of the solution on a larger cube controls the $L^{\infty}$ norm on a smaller cube, universally throughout the solution set.  A simple corollary of this result is that solutions bounded by $2$ in a large cube with sufficiently small energy must be bounded by a bit less than $2$ in the smaller cube, where the improvement is again universal. 

In the non-local setting, the $L^{2}$ to $L^{\infty}$ step works in essentially the same way as the local setting if one works in the whole space (instead of over finite cubes).  This is the content of Lemma \ref{L.L2toLinfWholeSpace} below.  The proof is included for completeness, and as a useful reference point for comparison with Proposition \ref{P.DiGeorgiIteration} that follows.

Unfortunately, Lemma \ref{L.L2toLinfWholeSpace} is not sufficient for implementing the grander oscillation reduction scheme of De Giorgi, since we need to track the solution on small cubes.  However, working locally in space introduces the long range effects of the fractional laplacian.  An insight of Caffarelli-Vasseur is that one can alternatively attack directly the usual corollary, but stated for the extension $\theta^{*}$ rather than $\theta$.  This is the content of Proposition \ref{P.DiGeorgiIteration} below.

\begin{Lem} \label{L.L2toLinfWholeSpace}
There exists a universal constant $C(d,\alpha)$ with the following property. For all $\theta \in \mathcal{S}(D,\alpha)$ with $\theta_{0} \in L^{2}(\R^{d})$ and every time $t>0$,
\begin{equation}
\|\theta \|_{L^{\infty}([t,\infty] \times \R^{d})} \leq C \|\theta_{0}\|_{L^{2}(\R^{d})}t^{-\frac{d}{4\alpha}}.
\end{equation}
\end{Lem}
\begin{Rem}
Note that $d$ should not be confused with $D$.  The constant $C$ in Lemma \ref{L.L2toLinfWholeSpace} is independent of the size of the critical drift $D$.
\end{Rem}
\begin{proof}
We will begin by reducing the Lemma to the following claim:
\begin{Claim}\label{Claim.ScalingRedLinf}
There exists a positive $\epsilon_{0}(\alpha,d)$ with the following property.  If $\theta \in \mathcal{S}(D,\alpha)$ and $\|\theta_{0}\|_{L^{2}(\R^{d})}^{2}\leq \epsilon_{0}$, then
$$
\|\theta\|_{L^{\infty}([1,\infty] \times \R^{d})} \leq 1.
$$
\end{Claim}
Assuming the claim, we now complete the proof. Let $\theta \in S(D,\alpha)$ and $\theta(0) \in L^{2}(\R^{d})$.  For each time $t>0$, we rescale and define $\tilde{\theta} \in \mathcal{S}(D,\alpha)$ via:
\begin{equation}\label{eq:RecalingLinfStep}
\tilde{\theta}(s,x)=t^{\frac{d}{4\alpha}}\epsilon_{0}^{\frac{1}{2}}\|\theta_{0}\|_{L^{2}(\R^{d})}^{-1}\theta(ts,t^{\frac{1}{2\alpha}}x).
\end{equation}
Applying the claim to $\tilde{\theta}$ and undoing the scaling gives:
\begin{equation}\label{eq:Einfty}
\|\theta\|_{L^{\infty}([t,\infty] \times \R^{d})} \leq \epsilon_{0}^{-\frac{1}{2}} \|\theta_{0}\|_{L^{2}(\R^{d})}t^{-\frac{d}{4\alpha}}.
\end{equation}
Setting $C=\epsilon_{0}^{-\frac{1}{2}}$ completes the proof of the Lemma.  Hence, it suffices to prove Claim \ref{Claim.ScalingRedLinf}.

\end{proof}

\begin{proof}[Proof of Claim~\ref{Claim.ScalingRedLinf}] 
We prove the claim by show that for $\epsilon_{0}>0$ sufficiently small, the inequality $\| \theta_{0}\|_{L^{2}(\R^{d})}\leq \epsilon_{0}$ implies
\begin{equation} \label{eq:ZeroTrunc}
||(\theta-1)_+||^{2}_{L^{\infty}([1,\infty);L^{2}(\R^{d})}=0
\end{equation}
Applying this observation to both $\theta$ and $-\theta$ completes the proof of the claim.

For each $k \in \N_{0}$,  we define:
\begin{itemize}
\item A time cut $T_{k}=1-2^{-k}$, a level $\lambda_{k}=1-2^{-k}$, and a truncation $\theta_{k}=(\theta - \lambda_{k})_{+}$. 
\item An energy $E_{k}=\|\theta_{k}\|^{2}_{L^{\infty}([T_{k},\infty);L^{2}(\R^{d})}+\|\theta_{k}\|^{2}_{L^{2}([T_{k},\infty];H^{\alpha}(\R^{d}))}$.
\end{itemize}
Our goal is to establish a nonlinear recursive relation between the energy at step $k$ and the energy at step $k-1$.  The first point is to observe that for all times $0 \leq s < t$, we know that
\begin{equation} \label{E.EnergyIneq}
\frac{1}{2}\|\theta_{k}(t)\|_{L^{2}(\R^{d})}^{2} + \|\theta_{k}\|^{2}_{L^{2}([s,t];\dot{H}^{\alpha}(\R^{d}))} \leq \frac{1}{2}\|\theta_{k}(s)\|_{L^{2}(\R^{d})}^{2}.
\end{equation}

The next point is to apply \eqref{E.EnergyIneq} in two different ways.  First we work with $k=s=0$ and maximize over $t>0$ to find.
\begin{equation} \label{E.IniitalSmallness}
\begin{split}
E_{0} &\leq \|\theta_{+}\|_{L^{\infty}([1/2,\infty);L^{2}(\R^{d}))}^{2} + \|\theta_{+}\|_{L^{2}([1/2,\infty);\dot{H}^{\alpha}(\R^{d}))}^{2} 
 \leq 2\|\theta_{+}(0)\|_{L^{2}(\R^{d})}^{2} \leq 2\epsilon_{0}.
\end{split}
\end{equation}
The second application of \eqref{E.EnergyIneq} is with an arbitrary $k \in \N$.  Fixing an $s \in [T_{k-1},T_{k}]$ then maximizing over $t \in [T_{k},\infty)$ yields 
$$
E_{k}\leq 2\|\theta_{k}(s)\|_{L^{2}(\R^{d})}^{2}\;\;\;\;\forall s\in[T_{k-1},T_{k}].
$$
Averaging this inequality over $s$ we find:
\begin{equation} \label{E.InitialEnergyBound}
E_{k}
\leq\min_{s\in[T_{k-1},T_{k}]} 2\|\theta_{k}(s)\|_{L^{2}(\R^{d})}^{2}
\leq \frac{1}{{T_{k}-T_{k-1}}}\int_{T_{k-1}}^{T_{k}} 2\|\theta_{k}(t)\|_{L^{2}(\R^{d})}^{2}\dt
\leq 2^{k+1}\|\theta_{k}\|^{2}_{L^{2}([T_{k-1},\infty) \times \R^{d})}.
\end{equation}
The next step is to interpolate and use a Sobolev embedding to find:
\begin{equation} \label{E.Interpolate}
\begin{split}
&\|\theta_{k-1}\|_{L^{2+\frac{4\alpha}{d}}\big ([T_{k-1},\infty) \times \R^{d}\big )} 
\leq \|\theta_{k-1}\|_{L^{\infty}  ([T_{k-1},\infty);L^{2}(\R^{d}))}^{\frac{2\alpha}{d+2\alpha}}
\|\theta_{k-1}\|_{L^{2} ([T_{k-1},\infty);L^{\frac{2d}{d-2\alpha}}(\R^{d}))}^{\frac{d}{d+2\alpha}}\\
& \leq E_{k-1}^{\frac{\alpha}{d+2\alpha}}\|\theta_{k-1}\|_{L^{2} ([T_{k-1},\infty);\dot{H}^{\alpha}(\R^{d}))}^{\frac{d}{d+2\alpha}}\leq E_{k-1}^{\frac{1}{2}}.
\end{split}
\end{equation}
Observe that if $\theta_{k}(t,x)>0$, then $\theta_{k-1}(t,x) > 2^{-k}$.  Hence, $1_{\{ \theta_{k}>0 \}}\leq 2^{\frac{4\alpha k}{d}}\theta_{k-1}^{\frac{4\alpha}{d}}$.  Combining this observation with inequalities \eqref{E.InitialEnergyBound} and \eqref{E.Interpolate}, we find that:
\begin{equation} \label{E.NonLinearIneq}
\begin{split}
E_{k} &\leq 2^{k+1}\|\theta_{k}1_{\{\theta_{k}>0\}}\|^{2}_{L^{2}([T_{k-1},\infty] \times \R^{d})}\\
& \leq 2^{(\frac{4\alpha}{d}+1)k+1} \|\theta_{k-1}\|_{L^{2+\frac{4\alpha}{d}}\big ([T_{k-1},\infty) \times \R^{d}\big )}^{2+\frac{4\alpha}{d}}\\
& \leq 2^{(\frac{4\alpha}{d}+1)k+1} E_{k-1}^{1+\frac{2\alpha}{d}}.
\end{split}
\end{equation}
Claim \ref{Claim.ScalingRedLinf} now follows from \eqref{E.IniitalSmallness} and \eqref{E.NonLinearIneq} for $\epsilon_{0}$ sufficiently small. In particular, one can prove there exists a universal $M>1$, such that $E^k\le M^{-k}$ for all $k \in \N$. Passing $k \to \infty$ yields \eqref{eq:ZeroTrunc}, for $\epsilon_0$ small enough.
\end{proof}

Now we move towards a variant of Lemma \ref{L.L2toLinfWholeSpace}, but working locally on cubes of the form $B_{r}=[-r,r]^{d}$.  The key starting point for the proof of Lemma \ref{L.L2toLinfWholeSpace} is the global energy inequality \eqref{E.EnergyIneq}, applied on truncations.  Naturally, one would like to prove a nonlinear inequality of the form \eqref{E.NonLinearIneq}, but replacing the role of the global energy inequality by the Cacciopoli type inequality \eqref{eq:NoCutoffCacciopoli}.  The main difficulty in implementing this idea within the general strategy of Lemma \ref{L.L2toLinfWholeSpace} is the contribution of $\theta^{*}$, which appears on the RHS of inequality \eqref{eq:NoCutoffCacciopoli}.

To understand the the proof of Proposition \ref{P.DiGeorgiIteration} below, a first observation (which can be checked with a small computation), is that the arguments in Lemma \ref{L.L2toLinfWholeSpace} would be robust if we could replace $\theta^{*}$ by $(\eta \theta)^{*}$ in \eqref{eq:NoCutoffCacciopoli}, for an appropriate cutoff $\eta$, at each step of the iteration.  This would effectively block out the non-local effects of the fractional laplacian.  However, a pointwise inequality of the form $\theta^{*} \leq (\eta \theta)^{*}$ cannot hold generically on a cube, even under the hypothesis $\theta^{*} \leq 2$, unless we add a barrier $b$ to the RHS.  Moreover, if we add a barrier $b_{k}$ at every step of the iteration, this could accumulate in the nonlinear inequality \eqref{E.NonLinearIneq} and kill the scheme.  An important observation of Caffarelli-Vasseur is the following: if the De Giorgi scheme is working, meaning the energy is decaying geometrically, and we shrink the $z$ domain (decreasing the influence of the barrier), then the truncation from a previous step will knock out the influence of the barrier at the present step.  This allows one to obtain the non-linear inequality in the next step and propagate the decay of energy.

The proof of Proposition \ref{P.DiGeorgiIteration} then splits roughly into three parts.  The first part is to proving Lemma \ref{L.PropOfSupport}, which shows that if the scheme is working, then we can replace the $k+1$ truncation $\theta_{k+1}^{*}$ by the extension of the $k$ truncation $(\eta_{k}\theta_{k})^{*}$ in the Cacciopoli inequality.  The second step is the proof that conditional on such a replacement, a nonlinear inequality for the energy holds.  This is the content of Lemma \ref{L.ConditionalNonlinearIneq}.  The final step is to inductively play these two Lemmas off of each other to complete the proof.


\begin{Prop}\label{P.DiGeorgiIteration}
There exist positive universal constants $\epsilon_{0}$ and $\lambda_{0}$ such that for all $\theta \in \mathcal{S}(D,\alpha)$, the following statement is true.   If $\theta^{*} \leq 2$ in $B_{3}^{*} \times [1,6]$ and
\begin{equation}
\|\theta_{+}^{*}\|^{2}_{L^{2}_{\omega_{\alpha}}([4,6] \times B_{2}^{*})} + \|\theta_{+}\|^{2}_{L^{\frac{2\alpha}{1-\alpha}}([4,6];L^{2}(B_{2}) )} \leq \epsilon_{0},
\end{equation}
then $\theta^{*} \leq 2-\lambda_{0}$ in $B_{1}^{*} \times [5,6]$.
\end{Prop}

We will begin by defining a sequence of energies $E_{k}$, in terms of two parameters $\lambda, \delta \in (0,1)$.  These constants will be chosen in Section \ref{SS.ProofofDiGeorgiIteration} at the end of the proof.  

For each $k \in \N$, we define:
\begin{itemize}
\item Radii $r_{k}=5/4 + 2^{-k}$, times $T_{k}=5-2^{-k}$, and levels $\lambda_{k}=\lambda(1+2^{-k})$.

\item Rectangles $B_{k}^{*}=B_{r_{k}} \times [0,\delta^{k}]$, $Q_{k}=[T_{k},6] \times B_{r_{k}}$, and $Q_{k}^{*}=B_{r_{k}}^{*} \times [T_{k},6]$.

\item Cutoffs $\sigma_{k} \in C^{\infty}_{c}(\R^{d} \times \R_{+})$ such that $1_{B_{k+1}^{*}} \leq \sigma_{k} \leq 1_{B_{k}^{*}}$ with $\eta_{k}(x)=\sigma_{k}(x,0)$ the trace of $\sigma_{k}$.

\item Truncations $\theta_{k}=(\theta-2+\lambda_{k})_{+}$ and truncations of the extension $\theta_{k}^{*}=(\theta^{*}-2+\lambda_{k})_{+}$.

\item Energies $E_{k}=\|\theta_{k}\|^{2}_{L^{\infty}([T_{k},6]; L^{2}(B_{k}))} + \|\nabla(\sigma_{k}\theta_{k}^{*})\|^{2}_{L^{2}_{\omega_{\alpha}}([T_{k},6];L^{2}(B_{k}^{*}))}$.
\end{itemize}


\subsection{Barriers} \label{SS.Barriers}
\begin{Lem} \label{L.BarrierInductive}
For each $k \in \N$, there exists a barrier $b_{k}: B_{k}^{*} \to \R$, a super-solution to
\begin{equation}\label{eq:equationforthebarrier}
\left\{
\begin{array}{rcl}
\Div(z^{1-2\alpha}\nabla b_k)&=0&\mbox{in \quad $B_{k}^{*}$}\\
b_k& = 2&\mbox{on \quad $\partial B_{k}^{*} \setminus  \{z=0\} \cup \{z=\delta^{k}\} $}\\
b_k&=0&\mbox{on \quad $ \{z=0\} \cup \{z=\delta^{k}\}$},
\end{array}\right.
\end{equation}
which satisfies, for universal constants $\overline{C},\;\mu > 0$, the inequality
\begin{equation} \label{E.BarrierDecay}
\sup_{(x,z) \in B_{k+1}^{*}} b_{k}(x,z) \leq \overline{C}e^{-\frac{\mu}{(2\delta)^{k}}}.
\end{equation}
\end{Lem}
\begin{proof}
We start by constructing a one dimensional barrier $b: \R_{+} \times [0,1] \to \R$, then use it to construct $b_{k}$. Namely, we build a super-solution $b(x,z)$ to
\begin{equation}
\left\{
\begin{array}{rcl}
\Div_{x,z}(z^{1-2\alpha}\nabla b)&\le 0&\mbox{in \quad $\R_{+} \times (0,1)$}\\
b&\ge 2&\mbox{on \quad $\{ x=0 \} $}\\
b&\ge0&\mbox{on \quad $ \{z=0\} \cup \{z=1\}$}.
\end{array}\right.
\end{equation}
In fact, by direct computations, one can check that
\begin{equation} \label{E.UnitScaleDecay}
b(x,z)=2\frac{e^{-\tilde{\mu} x}\cos(Az)}{\cos(A)}.
\end{equation}
is sufficient provided that $A \in (0,\pi/2)$ and $A^{2}\big ( (2\alpha-1)-\cos(A) \big)+\tilde{\mu}^{2} \leq 0$.  Next, define $b_{k}: B_{k}^{*} \to \R$
\begin{equation}
b_{k}(x,z)= \sum_{i=1}^{d}b\left(\frac{x_{i}+r_{k}}{\delta^{k}},\frac{z}{\delta^{k}}\right)+b\left(\frac{-x_{i}+r_{k}}{\delta^{k}},\frac{z}{\delta^{k}}\right).
\end{equation}
The scaling is has been chosen to ensure $b_{k}$ is a super solution to \eqref{eq:equationforthebarrier}. In view of \eqref{E.UnitScaleDecay}, we find that:
\begin{equation}
\sup_{(x,z) \in B_{k+1}^{*}} b_{k}(x,z)\leq \frac{4d}{\cos(A)} e^{-\tilde{\mu}\frac{(r_{k}-r_{k+1})}{\delta^{k}}}.
\end{equation}
Observing that $r_{k+1}-r_{k}=2^{-(k+1)}$, inequality \ref{E.BarrierDecay} follows if we set $\overline{C}=\frac{4d}{\cos(A)}$ and $\mu = \tilde{\mu}/2$.
\end{proof}
\subsection{Propagation of support} \label{SS.PropofSupp}
\begin{Lem}\label{L.PropOfSupport}
Assume that $\theta_{k}^{*}$ vanishes in the set $\big\{(t,x,z) \in Q_{k}^{*} \mid z=\delta^{k}\big\}$. In addition, suppose that the following two inequalities hold:
\begin{align} 
\label{H.TruncationKillsBarrier}\overline{C}e^{-\frac{\mu}{(2\delta)^{k}}} & \leq \lambda 2^{-k-2}.\\
\label{H.EnergyDecaysFast}E_{k}^{\frac{1}{2}} &\leq \lambda 2^{-k-2}\|P^{\alpha}_{\delta^{k+1}}\|_{L^{2}(\R^{d})}^{-1}.
\end{align}
Then
\begin{equation} \label{eq:keyInequality}
\theta_{k+1}^{*}1_{Q_{k+1}^{*}}\leq(\eta_{k}\theta_{k})^{*}.
\end{equation}
Moreover, $\theta_{k+1}^{*}$ vanishes in the set $\big\{(t,x,z) \in Q_{k+1}^{*} \mid z=\delta^{k+1}\big\}$.
\end{Lem}
\begin{proof}
Observe that the barrier $b_{k}$ has been constructed such that $\eta_{k}\theta_{k}^{*}$ is dominated by \linebreak $(\eta_{k}\theta_{k})^{*}+b_{k}$ in the region $Q_{k}^{*}$ by the maximum principle.  This follows from considering each portion of the boundary and using the support hypothesis for $\theta_{k}^{*}$.  Moreover, moving a bit further inside this region, the decay estimate \eqref{E.BarrierDecay} combined with \eqref{H.TruncationKillsBarrier} yield the following inequalities in $Q_{k+1}^{*}$:
\begin{equation} \label{E.MaxPrinciple}
\begin{split}
\eta_{k}\theta_{k}^{*} &\leq (\eta_{k}\theta_{k})^{*}+b_{k} \\
&\leq (\eta_{k}\theta_{k})^{*}+ \overline{C}e^{-\frac{\mu}{(2\delta)^{k}}}\\
&\leq (\eta_{k}\theta_{k})^{*}+ \lambda 2^{-k-2}.
\end{split}
\end{equation}
Now observe that in the region $Q_{k+1}^{*} \cap \{ \theta_{k+1}^{*}>0\}$, we have the identity $\theta_{k+1}^{*}=\eta_{k}\theta_{k}^{*}-\lambda 2^{-k-1}$.  Hence, multiplying both sides of \eqref{E.MaxPrinciple} above by $1_{Q_{k+1}^{*} \cap \{ \theta_{k+1}^{*}>0\} }$ and subtracting $\lambda 2^{-k-1}1_{Q_{k+1}^{*} \cap \{ \theta_{k+1}^{*}>0\} }$ yields the inequality:
\begin{equation}\label{eq:PropagationEst}
\begin{aligned}
\theta_{k+1}^{*}1_{Q_{k+1}^{*} \cap \{ \theta_{k+1}^{*}>0\}} &\leq (\eta_{k}\theta_{k})^{*}+ (\lambda 2^{-k-2}-\lambda 2^{-k-1})1_{Q_{k+1}^{*} \cap \{ \theta_{k+1}^{*}>0\}} \\
&\leq (\eta_{k}\theta_{k})^{*}.
\end{aligned}
\end{equation}
The inequality \eqref{eq:PropagationEst} now implies \eqref{eq:keyInequality} since the inequality holds trivially on the set $\{ \theta_{k+1}^{*}=0\}$ as $(\eta_{k}\theta_{k})^{*} \geq 0$.\\

\noindent Next we apply Young's inequality for convolutions, the definition of the energy $E_k$, and \eqref{H.EnergyDecaysFast} to obtain: 
\begin{align}
\sup_{(t,x) \in Q_{k}}\big\|(\eta_{k}\theta_{k})^{*}(t,x,\delta^{k+1})\big\|_{L^{\infty}(\R^{d})}
&\leq \|P^{\alpha}_{\delta^{k+1}}\|_{L^{2}(\R^{d})}\sup_{t \in [T_{k},4]}\big\|\eta_k\theta_k(t)\big\|_{L^2(\R^d)}\\
&\leq \|P^{\alpha}_{\delta^{k+1}}\|_{L^{2}(\R^{d})}E_{k}^{\frac{1}{2}} 
\leq \lambda 2^{-k-2}.
\end{align}
Combining this with \eqref{E.MaxPrinciple}, we can conclude that $\theta_{k+1}^{*}$ vanishes in the region \linebreak $\{(t,x,z) \in Q_{k+1}^{*} \mid z=\delta^{k+1}\}$.
\end{proof}
\subsection{Conditional non-linear inequality} \label{SS.CondNonlin}
\begin{Lem} \label{L.ConditionalNonlinearIneq}
There exist positive universal constants $\beta$ and $C$ such that for all $\theta \in \mathcal{S}(D,\alpha)$, the following statement holds true.
Let $k\geq 4$ and assume that 
\begin{equation} \label{E.ConditionNeededToDecay}
\theta_{k-1}^{*}1_{Q_{k-1}^{*}} \leq (\eta_{k-2}\theta_{k-2})^{*},
\end{equation}
then the following nonlinear inequality holds: 
\begin{equation}
E_{k} \leq C \left [ \frac{2^{(2\beta+4)k}}{\lambda^{2\beta}\delta^{2k}} \right ]E_{k-3}^{1+\beta}.
\end{equation}
\end{Lem}
\begin{proof}
The first step of the proof is to observe that there is a universal constant $C>0$ such that:
\begin{equation} \label{E.FirstEstforEk}
E_{k}\leq C2^{4k}\delta^{-2k} \bigg [ \|\theta_{k}^{*}\|^{2}_{L^{2}_{\omega_{\alpha}}(Q_{k-1}^{*})}+\|\theta_{k}\|^{2}_{L^{\frac{2\alpha}{1-\alpha}}([T_{k-1},6];L^{2}(B_{k-1}))} \bigg ].
\end{equation}
In order to establish \eqref{E.FirstEstforEk}, note first that:
\begin{equation}
\|\nabla(\sigma_{k}\theta_{k}^{*})\|^{2}_{L^{2}_{\omega_{\alpha}}(Q_{k}^{*})} \leq 2\|\nabla \sigma_{k}\|^{2}_{L^{\infty}(\R^{d} \times \R_{+})} \|\theta_{k}^{*}\|^{2}_{L^{2}_{\omega_{\alpha}}(Q_{k}^{*})} + 2\|\nabla\theta_{k}^{*}\|^{2}_{L^{2}_{\omega_{\alpha}}(Q_{k}^{*})}.
\end{equation}
Note that $\|\nabla \sigma_{k}\|^{2}_{L^{\infty}(\R^{d} \times \R_{+})} \leq C2^{2k}\delta^{-2k}$.  It now follows from the definition of $E_{k}$ that:
\begin{equation} \label{E.CutoffGrad}
E_{k} \leq C2^{2k}\delta^{-2k}\|\theta_{k}^{*}\|^{2}_{L^{2}_{\omega_{\alpha}}(Q_{k}^{*})}
+\|\theta_{k}\|^{2}_{L^{\infty}([T_{k},6];B_{k})}+2\|\nabla\theta_{k}^{*}\|^{2}_{L^{2}_{\omega_{\alpha}}(Q_{k}^{*})}.
\end{equation} 
\noindent Now we apply Proposition \ref{C.NoCutoffCacciopioli} with $r_{1}=r_{k}$, $\delta_{1}=\delta^{k}$ and $r_{2}=r_{k-1}$, $\delta_{2}=\delta^{k-1}$.  Working with times $t_{1} \in [T_{k-1},T_{k}]$ and $t_{2} \in [T_{k},6]$; maximizing in $t_{2}$ then averaging in $t_{1}$(as in the proof of Lemma \ref{L.L2toLinfWholeSpace}), we find that:
\begin{equation} \label{E.AppOfCorr}
\begin{split}
&\|\theta_{k}\|^{2}_{L^{\infty}([T_{k},6];B_{k})}+2\|\nabla\theta_{k}^{*}\|^{2}_{L^{2}_{\omega_{\alpha}}(Q_{k}^{*})}\\
&\leq C2^{k}\|\theta_{k}\|^{2}_{L^{2}(Q_{k-1})}+C2^{4k}\delta^{-2k}\bigg [ \|\theta_{k}^{*}\|^{2}_{L^{2}_{\omega_{\alpha}}(Q_{k-1}^{*})}+\|\theta_{k}\|^{2}_{L^{\frac{2\alpha}{1-\alpha}}([T_{k-1},6];L^{2}(B_{k-1}))} \bigg].
\end{split}
\end{equation}
Observe that Holder's inequality (in time) implies 
\begin{equation}\label{eq:TimeHolder}
\|\theta_{k}\|^{2}_{L^{2}(Q_{k-1})} \leq C_{\alpha}\|\theta_{k}\|^{2}_{L^{\frac{2\alpha}{1-\alpha}}([T_{k-1},6];L^{2}(B_{k-1}))}.
\end{equation}
Choosing an appropriate $C$, we may combine \eqref{E.CutoffGrad}, \eqref{E.AppOfCorr}, and \eqref{eq:TimeHolder} to obtain \eqref{E.FirstEstforEk}. \\

\noindent The second step of the proof is to deduce that:
\begin{equation}\label{E.ReductiontoNoStars}
\|\theta_{k}^{*}\|^{2}_{L_{\omega_{\alpha}}^{2}(Q_{k-1}^{*})}
+
\|\theta_{k}\|^{2}_{L^{\frac{2\alpha}{1-\alpha}}([T_{k-1},6];L^{2}(B_{k-1}))} 
\leq
 C \|\eta_{k-2}\theta_{k-2}\|^{2}_{L^{\frac{2\alpha}{1-\alpha}}([T_{k-2},6];L^{2}(B_{k-2}))}.
\end{equation}
To prove \eqref{E.ReductiontoNoStars}, we note that $\|\theta_{k}^{*}\|^{2}_{L_{\omega_{\alpha}}^{2}(Q_{k-1}^{*})} \leq \|\theta_{k-1}^{*}\|^{2}_{L_{\omega_{\alpha}}^{2}(Q_{k-1}^{*})}$. Hence, invoking our key assumption, inequality \eqref{E.ConditionNeededToDecay}, we find:
\begin{align}
\|\theta_{k}^{*}\|^{2}_{L_{\omega_{\alpha}}^{2}(Q_{k-1}^{*})}&\leq \|(\eta_{k-2}\theta_{k-2})^{*}\|^{2}_{L_{\omega_{\alpha}}^{2}(Q_{k-1}^{*})}\\
&\leq \left(\int_0^{\delta^{k-1}}z^{1-2\alpha}\|P_z^\alpha\|_{L^1(\R^d)}\dz\right)\|\eta_{k-2}\theta_{k-2}\|^{2}_{L^{2}(Q_{k-1})}\\
&=\frac{1}{2(1-\alpha)} \delta^{2(k-1)(1-\alpha)}\|\eta_{k-2}\theta_{k-2}\|^{2}_{L^{2}(Q_{k-1})}\\
&\leq C \|\eta_{k-2}\theta_{k-2}\|^{2}_{L^{\frac{2\alpha}{1-\alpha}}([T_{k-2},4];L^{2}(B_{k-2}))}.
\end{align} 
In the second inequality we used Young's inequality for the convolution and in the third inequality we used the fact that $\|P^{\alpha}_{z}\|_{L^{1}(\R^{d})}=1$ for all $z>0$. In the fourth inequality, we used that $\delta<1$.  The desired inequality \eqref{E.ReductiontoNoStars} now follows from the fact that 
\begin{equation}
\|\theta_{k}\|^{2}_{L^{\frac{2\alpha}{1-\alpha}}([T_{k-1},6];L^{2}(B_{k-1}))} \leq \|\eta_{k-2}\theta_{k-2}\|^{2}_{L^{\frac{2\alpha}{1-\alpha}}([T_{k-2},6];L^{2}(B_{k-2}))}.
\end{equation}
\noindent The third step in the proof is independent of the first two steps.  We aim to find an exponent $p=p(\alpha)>2$ such that
\begin{equation}\label{eq:raisinghomogeneity}
\|\eta_{k-3}\theta_{k-3}\|_{L^{\frac{4\alpha}{1-\alpha}}([T_{k-3},6];L^{p}(B_{k-3}))} \leq E_{k-3}^{\frac{1}{2}}.
\end{equation}
\noindent Observe that the variational form of the $H^{\alpha}(\R^{d})$ norm together with the fact that $\sigma_{k-3}$ is supported in $Q_{k-3}^{*}$ imply that:
\begin{equation}
\|\eta_{k-3}\theta_{k-3}\|_{L^{2}([T_{k-2},6];H^{\alpha}(\R^{d}))}^{2}
\leq \big\|\nabla (\sigma_{k-3}\theta_{k-3}^{*})\big\|_{L_{\omega_{\alpha}}^{2}(Q_{k-3}^{*})}^{2}
\leq E_{k-3}.
\end{equation}
Interpolating the spaces $\LL{\infty}{2}$ and $\LL{2}{\frac{2d}{d-2\alpha}}$ with  parameter $(1-\alpha)/2\alpha$ yields $\LL{\frac{4\alpha}{1-\alpha}}{p}$, where $p(\alpha)$ is given by the relation
\begin{equation}
\frac{1}{p}=\frac{1}{2}\bigg(\frac{1}{\alpha}-1\bigg)\bigg(\frac{1}{2}-\frac{\alpha}{d}\bigg)+\frac{1}{2}\bigg(3-\frac{1}{\alpha}\bigg)\frac{1}{2}.
\end{equation}
Note that since $\alpha \in (1/2,1)$ the interpolation parameter is in $(0,1)$, hence $p(\alpha)>2$.

Combining the interpolation with the Sobolev embedding theorem gives:
\begin{align}
&\|\eta_{k-3}\theta_{k-3}\|_{L^{\frac{4\alpha}{1-\alpha}}([T_{k-3},6];L^{p}(B_{k-3}))}\\ 
&\leq \|\eta_{k-3}\theta_{k-3}\|_{L^{\infty}([T_{k-3},6];L^{2}(B_{k-3}))}^{\frac{1}{2}(3-\frac{1}{\alpha})}\|\eta_{k-2}\theta_{k-2}\|_{L^{2}([T_{k-2},6];L^{\frac{2d}{d-2\alpha}}(B_{k-3}))}^{\frac{1}{2}(\frac{1}{\alpha}-1)} \\
&\leq E_{k-3}^{\frac{1}{4}(3-\frac{1}{\alpha})}\|\eta_{k-3}\theta_{k-3}\|_{L^{2}([T_{k-3},6];H^{\alpha}(\R^{d}))}^{\frac{1}{2}(\frac{1}{\alpha}-1)}\leq E_{k-3}^{\frac{1}{2}}.
\end{align}
\noindent This completes the proof of \eqref{eq:raisinghomogeneity}.\\  

\noindent We are now prepared to complete the proof of the Lemma. Due to the difference in the exponents of integrability in time and space, we define the exponent $q(\alpha)=\min\{4,p(\alpha)\}>2$ and observe that:
\begin{equation}
1_{ \{\eta_{k-2} \theta_{k-2}>0 \} } \leq \left(\frac{2^{k}}{4\lambda}\right)^{(\frac{q}{2}-1)}(\eta_{k-3}\theta_{k-3})^{(\frac{q}{2}-1)}.
\end{equation}
Using this inequality we obtain:
\begin{align}
\|\eta_{k-2}\theta_{k-2}\|^{2}_{L^{\frac{2\alpha}{1-\alpha}}([T_{k-2},6];L^{2}(B_{k-2}))}&=\|1_{ \{ \eta_{k-2}\theta_{k-2}>0 \} } \eta_{k-2}\theta_{k-2}\|^{2}_{L^{\frac{2\alpha}{1-\alpha}}([T_{k-2},6];L^{2}(B_{k-2}))}\\
&\leq \left(\int_{T_{k-2}}^6\left(\int_{B_{k-2}}\left|\left(\frac{2^{k}}{4\lambda}\right)^{(\frac{q}{2}-1)}(\eta_{k-3}\theta_{k-3})^{(\frac{q}{2}-1)} \eta_{k-2}\theta_{k-2}\right|^2\dx \right)^{\frac{\alpha}{1-\alpha}}\dt\right)^\frac{1-\alpha}{\alpha}\\
&\leq \left(\frac{2^{k}}{4\lambda}\right)^{q-2} \left(\int_{T_{k-2}}^6\left(\int_{B_{k-3}}|\eta_{k-3}\theta_{k-3}|^q\dx \right)^{\frac{\alpha}{1-\alpha}}\dt\right)^\frac{1-\alpha}{\alpha}\\
&=\left(\frac{2^{k}}{4\lambda}\right)^{(q-2)}\|\eta_{k-3}\theta_{k-3}\|^{q}_{L^{\frac{q\alpha}{1-\alpha}}([T_{k-3},6];L^{q}(B_{k-3}))}
\end{align}
Now, there are two cases: if $q=4$ then $q\leq p$ and we apply H\"{o}lder in space.  If $q=p$, then $p\alpha/(1-\alpha)\leq 4\alpha/(1-\alpha)$ and we apply H\"{o}lder in time.  In either case, \eqref{eq:raisinghomogeneity} yields:
\begin{equation}\label{eq:lastinequality}
\begin{aligned}
\left(\frac{2^{k}}{4\lambda}\right)^{(q-2)}\|\eta_{k-3}\theta_{k-3}\|^{q}_{L^{\frac{q\alpha}{1-\alpha}}([T_{k-3},6];L^{q}(B_{k-3}))}
&\leq C \left(\frac{2^{k}}{4\lambda}\right)^{(q-2)}\|\eta_{k-3}\theta_{k-3}\|^{q}_{L^{\frac{4\alpha}{1-\alpha}}([T_{k-3},6];L^{p(\alpha)}(B_{k-3}))}\\
&\leq C\left(\frac{2^{k}}{4\lambda}\right)^{q-2}E_{k-3}^{\frac{q}{2}}=C\left(\frac{2^{k}}{4\lambda}\right)^{q-2}E_{k-3}^{1+\beta_{2}},
\end{aligned}
\end{equation}
where $\beta=\frac{q}{2}-1$. Combining \eqref{E.FirstEstforEk}, \eqref{E.ReductiontoNoStars} and \eqref{eq:lastinequality}, we obtain the desired result.
\end{proof}
\subsection{Proof of Proposition} \label{SS.ProofofDiGeorgiIteration}
We begin by constructing an initial barrier $b_{0}$ which solves the following problem:
\begin{equation}
\left\{
\begin{array}{rcl}
\Div(z^{1-2\alpha}\nabla b_0)&=0&\mbox{in \quad $B_{3}^{*}$}\\
b_0&=2&\mbox{on \quad $\partial B_{3}^{*} \setminus  \{z=0\} $}\\
b_0&=0&\mbox{on \quad $ \{z=0\}$}.
\end{array}\right.
\end{equation}
By the strong maximum principle, we may select a universal $\lambda>0$ such that $b_{0} \leq 2-2\lambda$ in $B_{2}^{*}$.  The next step is to set the constant $\delta>0$, together with another constant $M>0$, which will be used to track the decay of $E_{k}$.  Namely, let us prove the following claim:
\begin{Claim}\label{claim:constants}
There exist positive universal constants $\delta$, $M$, and $k_{1}$ such that the following inequalities hold for all $k\in\N$
\begin{align}
\label{H.TruncationKillsBarrier2}\overline{C}e^{-\frac{\mu}{(2\delta)^{k}}} &\leq \lambda 2^{-k-2},\\
\label{H.EnergyDecaysFast2}M^{-\frac{k}{2}} &\leq \lambda 2^{-k-2}\|P^{\alpha}_{\delta^{k+1}}\|_{L^{2}(\R^{d})}^{-1},
\end{align}
and for all $k \geq k_{1}$
\begin{align}
\label{H.SchemeIsWorking} C\left [ \frac{2^{(2\beta+4)k}}{\lambda^{2\beta}\delta^{2k}} \right ]M^{(3-k)(1+\beta)}  &\leq M^{-k},
\end{align}
where $\overline{C}$, $\mu$, $\beta$ and $C$ are universal constants from Lemma~\ref{L.BarrierInductive} and Lemma~\ref{L.ConditionalNonlinearIneq}.
\end{Claim}
\begin{proof}[Proof of Claim~\ref{claim:constants}]
First choose $\delta>0$ sufficiently small to make the first inequality hold for all $k\in\N$ and fix it. Next choose $M$ sufficiently large to make the next two inequalities hold, considering that $\|P^{\alpha}_{\delta^{k+1}}\|_{L^{2}(\R^{d})}^{-1}=\delta^{(k+1)d/2}\|P^{\alpha}_1\|_{L^{2}(\R^{d})}^{-1}$. The size of $k_1$ is related to the size of $\beta$.
\end{proof}

We are now prepared to prove the following:
\begin{Claim}\label{claim:energydecay}
There exists a positive universal $\epsilon_{0}$ such that for all $k$, the energy decays as follows: $E_{k} \leq M^{-k}$.
\end{Claim}
\begin{proof}
The first step is to choose $\epsilon_{0}$ sufficiently small to ensure that $E_{k} \leq M^{-k}$ for all $k\le k_1$.  Using Corollary~\ref{C.NoCutoffCacciopioli} we argue exactly as in the proof of Lemma \ref{L.ConditionalNonlinearIneq} to find a universal constant $C$ such that:
\begin{align}
E_{k}&\leq C2^{4k}\delta^{-2k } \bigg [ \|\theta_{k}^{*}\|^{2}_{L_{\omega_{\alpha}}^{2}(Q_{k-1}^{*})}+\|\theta_{k}\|^{2}_{L^{\frac{2\alpha}{1-\alpha}}([T_{k-1},6];L^{2}(B_{k-1}))} \bigg ] \\
&\leq C2^{4k}\delta^{-2k } \bigg [ \|\theta_{k}^{*}\|^{2}_{L_{\omega_{\alpha}}^{2}(Q_{2}^{*})}+\|\theta_{k}\|^{2}_{L^{\frac{2\alpha}{1-\alpha}}([T_{k-1},6];L^{2}(B_{2}))} \bigg ]\\
&\leq C2^{4k}\delta^{-2k }\epsilon_{0}, 
\end{align}
for all $k\in \N$. Choosing $\epsilon_{0}$ sufficiently small (depending only on the universal constants $\delta$ and $M$ set above), gives the desired initial energy decay $E_{k} \leq M^{-k}$ for all $k\le k_1$.  Next observe that the maximum principle on $B_3^*$ implies that in $Q_{1}^{*}$
\begin{equation}
\eta_{1}\theta^{*} \leq (\eta_{1}\theta)^{*}+b_{0}.
\end{equation}
By Young's inequality for convolutions, it follows that:
\begin{align}
\sup_{t\in [T_{1},4]}\big\|\eta_{1}\theta^*(t,.,\delta)\big\|^{2}_{L^{\infty}(\R^{d})} 
\leq \|P^{\alpha}_{\delta}\|_{L^{2}(\R^{d})}E_{1}^{\frac{1}{2}}.
\end{align}
By taking $\epsilon_0$ small enough, we can make this last term less than $\lambda/2$, so that combining this with $b_{0}\leq 2-2\lambda$ gives $\theta^{*}\leq 2-3\lambda/2=2-\lambda_1$ in the region $\{(t,x,z) \in Q_{1}^{*} \mid z=\delta\}$.

Now that the inductive hypotheses has been set, we may assume for the purpose of induction that $E_{j} \leq M^{-j}$ for all $j\leq k-1$, where $k-1\ge k_1$.  We will use Lemmas \ref{L.PropOfSupport} and \ref{L.ConditionalNonlinearIneq} together with Claim~\ref{claim:constants} to prove that $E_{k} \leq M^{-k}$.

Indeed, the inductive hypothesis together with \eqref{H.TruncationKillsBarrier2} and \eqref{H.EnergyDecaysFast2} allows us to repeatedly apply Lemma \ref{L.PropOfSupport} for each $j \leq k-2$, and deduce that  $\theta_{j}^{*}$ vanishes in the region $\{(t,x,z) \in Q_{j}^{*} \mid z=\delta^{j}\}$.  In addition, this yields the inequality
\begin{equation}
\theta_{k-1}^{*}1_{Q_{k-1}^{*}}\leq(\eta_{k-2}\theta_{k-2})^{*}.
\end{equation}
This allows us to apply Lemma \ref{L.ConditionalNonlinearIneq} to find:
\begin{equation}
E_{k} \leq  C\left [ \frac{2^{(2\beta+4)k}}{\lambda^{2\beta}\delta^{2k}} \right ]E_{k-3}^{1+\beta}
\leq C\left [ \frac{2^{(2\beta+4)k}}{\lambda^{2\beta}\delta^{2k}} \right ]M^{(3-k)(1+\beta)}
\leq M^{-k},
\end{equation}
where we have used the inductive hypothesis and \eqref{H.SchemeIsWorking}.
\end{proof}
\begin{proof}[Proof of Proposition~\ref{P.DiGeorgiIteration}]
By the Claim~\ref{claim:energydecay}, we deduce that:
\begin{equation}
\lim_{k \to \infty}E_{k}=0.
\end{equation}
In particular, this yields $\theta \leq 2-\lambda$ in $[5,6] \times B_{\frac{5}{4}}$.  The last step of the proof is to construct one final barrier $b_{1}$ which satisfies:
\begin{equation}
\left\{
\begin{array}{rcl}
\Div(z^{1-2\alpha}\nabla b_1)&=0&\mbox{in \quad $B_{\frac{5}{4}}^{*}$}\\
b_1&=2&\mbox{on \quad $\partial B_{\frac{5}{4}}^{*} \setminus  \{z=0\} $}\\
b_1&=2-\lambda &\mbox{on \quad $ \{z=0\}$}
\end{array}\right.
\end{equation}
By the strong maximum principle, there exists a $\lambda_{0}$ such that $\theta^{*} \leq 2-\lambda_{0}$ in $[5,6] \times B_{1}^{*}$.  This completes the proof.
\end{proof}

\section{Isoperimetric Inequality}\label{S.Isoperim}
Ultimately, we would like to apply Proposition \ref{P.DiGeorgiIteration} to obtain an oscillation reduction result.  However, we will need to verify the smallness constraint in order for the Proposition to be useful.  For $\theta$ such that $\theta^{*} \leq 2$, one way to check the smallness would be to first show that the measure of the transition set, where $0<\theta^{*}<1$, is small, then attack separately the region where $1 \leq \theta^{*}\leq 2$.  The next lemma shows that as long as $\theta^{*}$ is negative for a sufficient proportion of space/time, the second step is redundant.  Indeed, if the first step held but the second step didn't, this would indicate a jump discontinuity.  The classical isoperimetric inequality quantifies the sense in which the dissipative bounds rule out this possibility.  Unfortunately, we have no control on $\partial_{t}\theta^{*}$, so making this argument rigorous requires a careful analysis.  This is the content of Proposition \ref{lem:Isoperimetric}.
\begin{Prop}\label{lem:Isoperimetric}
For all $\epsilon_{1}>0$, there exist a $\delta_{1}>0$ with the following property.  For all $\theta \in \mathcal{S}(D,\alpha)$ with $\theta^{*} \leq 2$ in $[1,6] \times B_{3}^{*}$ and
$
\left |\left\{ \theta^{*}\leq 0 \right\} \cap [1,6] \times B_{3}^{*}\right| \geq (1/2)\left |[1,6] \times B_{3}^*\right|,
$
if
\begin{equation}\label{eq:smallnessInBetween}
\left|\left\{  0<\theta^{*}< 1 \right\} \cap [1,6] \times B_{3}^{*} \right|\leq\delta_{1},
\end{equation}
then
\begin{equation}\label{eq:smallnessAboveOne}
\left |\left\{\theta^{*}\geq 1\right\} \cap [4,6] \times B_{2}^{*} \right |\leq \epsilon_{1}.
\end{equation}
\end{Prop}
\begin{proof}First observe that if we prove the Lemma for sufficiently small $\epsilon_{1}$, the general case easily follows.  We will take $\epsilon_{1}\ll1$ throughout, re-adjusting the upper bound a finite number of times in the course of the argument.

For each $t \in [1,6]$, we will study the measure of the sets
\begin{align}
	&\mathcal{A}(t) = \big\{ (x,z) \in B_{r}^{*} : \theta^{*}(t,x,z) \leq 0 \big\}, \\	
	&\mathcal{B}(t) = \big\{ (x,z) \in B_{r}^{*} : \theta^{*}(t,x,z) \geq 1 \big\}, \\
	&\mathcal{C}(t) = \big\{ (x,z) \in B_{r}^{*} : 0<\theta^{*}(t,x,z) < 1 \big\},
\end{align}
where $r \in (2,3)$ is a constant (depending only on the dimension $d$) chosen below.  In fact, we will mostly restrict our attention to $t$ in a set $I$ of ``good'' times defined by:
\begin{align}
I=\left \{t \in [1,6]:  \int_{B_{r}^{*}}|\nabla \theta_{+}^{*}(t)|^{2}\dx \dz \leq K, \quad |C(t)|\leq \delta_{2} \right\},
\end{align}
where $K$ is a large parameter and $\delta_{2}$ is a small parameter, both to be chosen depending on $\epsilon_{1}$.  

The general strategy of proof is to choose the parameters $r,\;K,\;\delta_{1},\;\delta_{2}$ to ensure that \eqref{eq:smallnessInBetween} implies the following hold simultaneously:
\begin{enumerate}
\item \label{item:GoodTimesMeasure} Most times in $[1,6]$ are good: $\left |I^{c} \cap [1,6] \right| \leq \epsilon_{1}/2$.
\item \label{item:KeyConc} At each good time $t \in I \cap [4,6]$, it holds that $|\mathcal{B}(t)| \leq \epsilon_{1}/4$.
\end{enumerate}  
Combining the two easily yields \eqref{eq:smallnessAboveOne}.  The proof of \ref{item:GoodTimesMeasure} is a simple consequence of the Chebyshev inequality.  The proof of \ref{item:KeyConc} is subtle.   Let us start with the first point and prove the following.
\begin{Claim}\label{claim:smallTimes}
There is a universal $C_{0}>0$ such that if $K\leq C_{0}\epsilon^{-1}_{1}$ and $\delta_{1} \leq (1/4) \delta_{2}\epsilon_{1}$, then
\begin{equation}\label{eq:SmallTimeMeasure}
\left |I^{c} \cap [1,6] \right| \leq (1/2)\epsilon_{1}.
\end{equation}
\end{Claim}  
\begin{proof}
Since $r<3$ is and $\theta^{*}\leq 2$ in $[1,6] \times B_{3}^{*}$, the Cacciopoli inequality, Proposition \ref{C.NoCutoffCacciopioli}, yields a universal $C>0$ such that
\begin{equation}
\int_{1}^{6}\int_{B_{r}^{*}} z^{1-2\alpha}|\nabla \theta_{+}^{*}|^{2}\dx\dz\dt \leq C.
\end{equation}
Note that for each $t$ fixed we have the trivial inequality 
\begin{equation}
\begin{aligned}
	 \int_{B_{r}^{*}} |\nabla \theta_{+}^{*}(t)|^{2}\dx \dz \leq r^{2\alpha-1}\int_{B_{r}^{*}} z^{1-2\alpha}|\nabla \theta_{+}^{*}(t)|^{2}\;\dee x\dee z,
\end{aligned}
\end{equation}
so an application of the Chebyshev inequality yields:
\begin{equation}\label{eq:sizeK}
    \left|\left\{t \in [1,6] : \| \nabla\theta^{*}_{+}(t)\|_{L^{2}(B_{r}^{*})}^{2}\geq K\right\}\right|
    \leq Cr^{2\alpha-1}K^{-1}.
\end{equation}
Applying Chebyshev's inequality once more and recalling that $\left|\left\{  0<\theta^{*}< 1 \right\} \cap [1,6] \times B_{3}^{*} \right|\leq\delta_{1}$,
\begin{equation}\label{eq:sizeC}
\begin{aligned}
\big|\big\{t \in [1,6]  : |\mathcal{C}(t)| \geq \delta_{2} \big\}\big|
&\leq \delta_{2}^{-1}\left | \left\{ 0<\theta^{*}<1\right \} \cap B_{r}^{*} \times [1,6] \right |\leq \delta_{2}^{-1}\delta_{1},
\end{aligned}
\end{equation}
since $B_{r}^{*}\subset B_{3}^{*}$.  Combining \eqref{eq:sizeK} and \eqref{eq:sizeC} and setting $C_{0}=(C/4)r^{1-2\alpha}$ completes the proof.
\end{proof}
In accordance with Claim \ref{claim:smallTimes}, we now fix the constants $K$ and $\delta_{1}$ as $K=C_{0}\epsilon^{-1}_{1}$, $\delta_{1}=(1/4) \delta_{2}\epsilon_{1}$, and let $\delta_{2}=\epsilon_{1}^{p}$ for some universal exponent $p$ to be chosen below.

Now we proceed to the more difficult part of the proof, point 2.  A key tool is the isoperimetric inequality (see Appendix A in \cite{caffarelli2010drift}), which implies that for each $t \in I$
\begin{equation} \label{eq:GeneralIsoperim}
|\mathcal{A}(t)||\mathcal{B}(t)| \leq |\mathcal{C}(t)|^{\frac{1}{2}}K^{\frac{1}{2}} \leq (\delta_{2}K)^{\frac{1}{2}}\leq (C_{0}C_{1})^{\frac{1}{2}}\epsilon_{1}^{\frac{1}{2}(p-1)}.
\end{equation}
To use \eqref{eq:GeneralIsoperim} to establish the second point, we will need lower bounds for $|\mathcal{A}(t)|$.  Our main step is to prove the following:

\begin{Claim}\label{claim:propInTime}
There exists a universal constant $\kappa$ (independent of $\epsilon_{1}$) with the following property $:$ if $p>13$ and $|\mathcal{A}(t_{0})| \geq 1/8$ for some $t_{0}$, then $|\mathcal{B}(t)| \leq \epsilon_{1}/4$ for all $t \in [t_{0},t_{0}+\kappa] \cap I$.
\end{Claim}
\begin{proof}
The proof of the claim has three steps.  The first step is to show that:
\begin{equation} \label{eq:Transfer_to_Theta}
\int_{B_{r}}\theta_{+}^{2}(t_{0})\dx \leq (1/4)|B_{r}|.
\end{equation}
To establish \eqref{eq:Transfer_to_Theta}, first consider $\theta^{*}_{+}$.  By \eqref{eq:GeneralIsoperim} and the assumption $\theta^{*} \leq 2$,
\begin{align}\label{eq:controlforcaccioppolinonlocalpart}
\int_{B_{r}^{*}} \big(\theta_{+}^{*}(t_0)\big)^{2}\;\dee x\dee z \leq 4 \big ( |\mathcal{B}(t_0)|+|\mathcal{C}(t_0)| \big ) \leq 4 \left ( (\delta_{2}K)^{\frac{1}{2}}+\delta_{2} \right )\leq \epsilon_{1}^{\frac{1}{4}(p-1)}.
\end{align}
To connect $\theta_{+}$ to $\theta^{*}_{+}$, apply the Fundamental Theorem of Calculus for each $z \in [0,r]$ to find
\begin{equation}\label{eq:FTC}
	\int_{B_{r}}\theta_{+}^{2}(t_{0},x)\;\dee x = \int_{B_{r}}(\theta_{+}^{*})^{2}(t_{0},x,z)\;\dee x-2\int_{0}^{z}\int_{B_{r}}\theta_{+}^{*}(t_{0},x,z)\partial_{z}\theta^{*}_{+}(t_{0},x,z)\;\dee x\dee z.
\end{equation}
Averaging both sides over $z \in [0,r]$ and using H\"older's inequality gives
\begin{equation}\label{eq:controlforcaccioppolilocalpart}
\begin{aligned}
	\int_{B_{r}}\theta_{+}^{2}(t_0)\dee x 
    & \leq r^{-1} \int_{B_{r}^{*}} (\theta_{+}^{*})^{2}(t_0)\dee x\dee z + 2 \int_{B_{r}^{*}}|\theta_{+}^{*}(t_0)|\partial_{z}\theta^{*}_{+}(t_0)|\dee x\dee z \\
	& \leq r^{-1}\int_{B_{r}^{*}} (\theta_{+}^{*})^{2}(t_0)\dee x\dee z + 2 \left( \int_{B_{r}^{*}}|\theta_{+}^{*}(t_0)|^{2}\dee x\dee z\right)^{\frac{1}{2}}  \left(\int_{B_{r}^{*}}|\nabla\theta_{+}^{*}(t_0)|^{2}\dee x\dee z \right)^{\frac{1}{2}}\\
    &\leq r^{-1} \epsilon_{1}^{\frac{1}{4}(p-1)} +  \epsilon_{1}^{\frac{1}{8}(p-1)}K^{\frac{1}{2}} 
     \leq \epsilon_{1}^{ (1/8)(p-1)-(1/2)}.
\end{aligned}
\end{equation}
The claim now follows as long as $p>5$.  

The second step is to choose $\kappa$ and check that for all $t \in [t_{0},t_{0}+\kappa]$,
\begin{equation} \label{eq:Propgate_in_Time}
\int_{B_{r}}\theta_{+}^{2}(t)\dx \leq (3/4)|B_{r}|.
\end{equation}
To this end, we apply the Caccioppoli estimate, Corollary~\ref{C.NoCutoffCacciopioli}, with the $(x,z)$ rectangles $B_{r} \times [0,1]$ and $B_{R} \times [0,2]$ on the time interval $[t_{0},t]$ for each $t<t_{0}+\kappa$.  This yields
\begin{equation}\label{eq:auxiliarysizeoftheL2norm}
\begin{aligned}
\int_{B_{r}}\theta_{+}^{2}(t)\;\dee x 
&\leq \int_{B_{r}} \theta_{+}^{2}(t_{0})\dee x 
+ \int_{B_{R} \setminus B_{r}} \theta_{+}^{2}(t_{0})\dee x 
+ C(R-r)^{-4}\big ( (t-t_{0})+(t-t_{0})^{\frac{1}{\alpha}-1} \big) \\
&\leq (1/4)|B_{r}|
+ 4|B_{R} \setminus B_{r}| 
+ C(R-r)^{-4}\kappa^{\frac{1}{\alpha}-1}. \\
\end{aligned}
\end{equation}
We may now choose $R$ close enough to $r$, then $\kappa$ sufficiently small (both independently of $\epsilon_{1}$) to deduce the desired inequality \eqref{eq:Propgate_in_Time}.

The final step is to complete the proof of the claim.  Towards this end, use again the connection \eqref{eq:FTC} of $\theta_{+}$ to $\theta^{*}_{+}$, so that for $z \in [0,r]$
\begin{align}
\int_{B_{r}} (\theta^{*}_{+})^{2}(t,x,z)\;\dee x \leq \int_{B_{r}}\theta_{+}^{2}(t,x)\dx + C\sqrt{z} K^\frac{1}{2} \leq (3/4)|B_{r}|+C\sqrt{z}K^{\frac{1}{2}}.
\end{align}
By Chebyshev's inequality,
$
|\{x \in B_{r} : \theta^{*}(t,x,z) \geq 1 \}| \leq (3/4)|B_{r}|+C\sqrt{z}K^{\frac{1}{2}}.
$
Hence we have
\begin{equation}\label{eq:sizeofthesetbiggerthanone}
\left |\left \{\theta^{*}(t) \geq 1 \right\} \cap B_{r} \times [0,z] \right| \leq z \left ((3/4)|B_{r}|+C(zK)^{\frac{1}{2}} \right ).
\end{equation}
Now, for $z\leq r$, we may bound $|\mathcal{A}(t)|$ from below via
\begin{equation}
\begin{aligned}
|\mathcal{A}(t)| &\geq \left |\left \{\theta^{*}(t) \leq 0 \right\} \cap B_{r} \times [0,z] \right | \\
&= |B_{r}|z -  \left |\left \{\theta^{*}(t) \geq 1 \right\} \cap B_{r} \times [0,z] \right|- \left| \left\{ 0<\theta^{*}(t)< 1 \right\} \cap B_{r} \times [0,z]\right| \\
& \geq |B_{r}|z-z\left ( (3/4) |B_{r}|+C(zK)^{\frac{1}{2}} \right )-\delta_{2}
\geq (1/4)|B_{r}|z-Cz^{\frac{3}{2}}/\epsilon_{1}^{\frac{1}{2}}-\epsilon_{1}^{p}.\\
\end{aligned}
\end{equation}
Choosing $z=\epsilon_{1}^{2}$, we find that 
$
|\mathcal{A}(t)|\geq C\epsilon_{1}^{2}.
$
Applying the Isoperimetric inequality we find that $|\mathcal{B}(t)| \leq C\epsilon_{1}^{\frac{1}{4}(p-1)-2}$.  Hence, for $p>13$ we find $|\mathcal{B}(t)| \leq \epsilon_{1}/4$ as desired.
\end{proof}
In view of Claim \ref{claim:propInTime}, we set $p=14$.  Finally, we choose $r$ in a way that yields a time $t_{0} \in [1,4)$ with $|\mathcal{A}(t_{0})| \geq 1/8$.  If no such time exists, then
\begin{equation}
\left|\left\{\theta^{*}>0\right\} \cap B_{3}^{*} \times [1,6]\right| \geq 3|B_{r}|^{*}-(3/8).
\end{equation}  
This gives a contradiction with our negativity hypothesis provided that $3|B_{r}|^{*}-(3/8)>(5/2)|B_{3}^{*}|$.  As the inequality holds strictly as $r \to 3$, we may choose an $r<3$ depending on the dimension $d$.  

Appealing to Claim \ref{claim:propInTime}, for all $t \in [t_{0},t_{0}+\kappa] \cap I$, we have $|\mathcal{B}(t)|\leq \epsilon_{1}/4$.  Moreover, in view of Claim \ref{claim:smallTimes}, we may ensure $|I^{c}|<\frac{\kappa}{2}$ by choosing $\epsilon_{1}<4\kappa$.  Thus, the set $I \cap [t_{0}+\kappa/2,t_{0}+\kappa]$ is non-empty and must contain a time $t_{1}$.  Hence, it follows that
\begin{equation}
|\mathcal{A}(t_{1})| \geq |B_{r}^{*}|-|\mathcal{B}(t_{1})|-|\mathcal{C}(t_{1})| \geq 1/8.
\end{equation}  
Applying the claim again yields $|\mathcal{B}(t)|\leq \epsilon_{1}/4$ for all $t \in [t_{1},t_{1}+\kappa] \cap I$. Repeating this argument finitely many times gives the second point, completing the proof.
\end{proof}
Finally, we use Proposition \ref{lem:Isoperimetric} to show that a large amount of negativity of $\theta^{*}$ combined with a small transition set is sufficient to decrease from $2$ to $2-\lambda_{0}$, where $\lambda_{0}$ comes from the Proposition \ref{P.DiGeorgiIteration}.
\begin{Cor}\label{Cor_Final_Isoperim}
There exists $\delta>0$ with the following property.  For all $\theta \in \mathcal{S}(D,\alpha)$ such that $\theta^{*}\leq 2$ and 
$\left|\left\{\theta^{*}\leq 0\right\} \cap B_{3}^{*} \times [1,6]\right | \geq (1/2)\left|B_{3}^{*} \times [1,6]\right|,$
the inequality 
\begin{equation}
\left|\left\{0<\theta^{*}<1\right\} \cap B_{3}^{*} \times [1,6]\right|<\delta
\end{equation}
implies that $\theta^{*}\leq 2-\lambda_{0}$ on $B_{1}^{*} \times [5,6]$.
\end{Cor}
\begin{proof}
Begin by applying the Cacciopioli inequality to deduce
\begin{equation}
\|\theta_{+}\|_{L^{\infty}([4,6];L^{2}(B_{2}))}+\int_{4}^{6}\int_{B_{2}^{*}}|\nabla\theta^{*}_{+}|^2\dx\dz\dt \leq C
\end{equation}
for some universal $C$ (with a simple argument to remove the weight $z^{1-2\alpha}$).  Next apply Lemma \ref{lem:Isoperimetric} with $\epsilon_{1}$ to be chosen.  If $\delta < \delta_{1}$ then
\begin{equation}
|\{\theta^{*}\geq 1\} \cap B_{2}^{*} \times [4,6]| \leq \epsilon_{1}.
\end{equation}
If, in addition $\delta<\epsilon_{1}$, we deduce further that
\begin{equation}
|\{\theta^{*}> 0\} \cap B_{2}^{*} \times [4,6]| \leq 2\epsilon_{1}.
\end{equation}
Applying H\"{o}lder's inequality with exponent $\alpha/(2\alpha-1) > 1$,
\begin{equation}\label{eq:nonlocalpartofisoperimetric}
\int_{4}^{6}\int_{B_{2}^{*}}z^{1-2\alpha}\left(\theta^{*}_{+}\right)^{2}\dz\dx\dt \leq C\epsilon_{1}^{\frac{1}{\alpha}-1}.
\end{equation}
Applying the Fundamental Theorem of Calculus,
\begin{align}
\int_{4}^{6}\int_{B_{2}}\theta_{+}^{2}\dx\dt &\leq \int_{4}^{6}\int_{B_{2}^{*}}\left(\theta^{*}_{+}\right)^{2}\dx\dz\dt +2\left (\int_{4}^{6}\int_{B_{2}^{*}}\left(\theta^{*}_{+}\right)^{2}\dx\dz\dt \right)^{\frac{1}{2}} \left(\int_{4}^{6}\int_{B_{2}^{*}}|\nabla\theta^{*}_{+}|^2\dx\dz\dt \right )^{\frac{1}{2}}\\ 
&\leq C \epsilon_{1}^{\frac{1}{2}(\frac{1}{\alpha}-1)}.
\end{align}
Interpolating in time with parameter $(1/\alpha)-1$ gives
\begin{equation}\label{eq:localpartofisoperimetric}
\begin{aligned}
\|\theta_{+}\|_{L^{\frac{2\alpha}{1-\alpha}}([4,6];L^{2}(B_{2}))}^{2}&\leq \|\theta_{+}\|_{L^{2}([4,6]\times B_{2})}^{\frac{2}{\alpha}-2}\|\theta_{+}\|_{L^{\infty}([4,6];L^{2}(B_{2}))}^{4-\frac{2}{\alpha}}\\
&\leq C\epsilon_{1}^{(\frac{1}{\alpha}-1 )^{2}}.
\end{aligned}
\end{equation}
Combining these observations, we find that
\begin{equation}
\|\theta_{+}^{*}\|^{2}_{L^{2}_{\omega_{\alpha}}([4,6] \times B_{2}^{*})} + \|\theta_{+}\|^{2}_{L^{\frac{2\alpha}{1-\alpha}}([4,6];L^{2}(B_{2}) )} \leq C\epsilon_{1}^{(\frac{1}{\alpha}-1 )^{2}}.
\end{equation}
Hence, setting $C\epsilon_{1}^{(\frac{1}{\alpha}-1 )^{2}}=\epsilon_{0}$ and $\delta<\delta_{1} \wedge \epsilon_{1}$, we may apply Proposition \ref{P.DiGeorgiIteration} and deduce that $\theta^{*} \leq 2-\lambda_{0}$ in $B_{1}^{*} \times [5,6]$.
\end{proof}


\section{Proof of Main Theorem}\label{S.Proof}
The next step is to remove the smallness hypothesis on the transition set required in Corollary \ref{Cor_Final_Isoperim}.  This follows from setting up the classical dichotomy from the local elliptic theory.  If the smallness of the transition set fails to hold, we can keep resizing the solution until it eventually does, or else we exhuast too much measure in space/time.  Moreover, to obtain a true oscillation reduction result, we need to remove the hypothesis $\theta^{*} \leq 2$ by another resizing argument.  The details are carried out for completeness in Proposition \ref{cor:reductionofoscillation} below.
\begin{Prop}\label{cor:reductionofoscillation}
There exist a positive universal constant $\mu<1$ such that for all $\theta \in \mathcal{S}(D,\alpha)$,
\begin{equation}\label{eq:decreaseofoscillations}
\underset{[5,6] \times B_{1}^{*}}{\text{osc}}\theta^{*}\le \mu \underset{[1,6] \times B_{3}^{*}}{\text{osc}} \theta^{*}.
\end{equation}
\end{Prop}
\begin{proof}
The first step of the proof is to remove the smallness constraint from Corollary \ref{Cor_Final_Isoperim} by establishing:
\begin{Claim}\label{eq:UnconditionDecreaseofLinf}
There exists a universal $\lambda>0$ such that for all $\theta \in \mathcal{S}(D,\alpha)$, the following statement is true.   If $\theta^{*} \leq 2$ in $B_{3}^{*}\times[1,6]$ and
$
|\{ \theta^{*} \leq 0 \} \cap B_{3}^{*}\times[1,6] | \geq (1/2)|B_{3}^{*}\times[1,6]|,
$
then $\theta^{*} \leq 2-\lambda$ in $B_{1}^{*} \times [5,6]$.
\end{Claim}
\begin{proof}
For $k \in \N$, define recursively $\theta_k=2(\theta_{k-1}-1), \, \mbox{with }\;\theta_0=\theta$.  Equivalently, 
\begin{equation}
\theta_{k}=2^k(\theta-2)+2.
\end{equation}  
Next we define a universal constant $K_{1}$ as the floor of $5\delta^{-1}|B_{3}^{*}|$, where $\delta$ was defined in Corollary \ref{Cor_Final_Isoperim}. Set $\lambda=\lambda_{0} 2^{-K_{1}}$. 

If it happens that for all $0\leq k\leq K_{1}$,
\begin{equation}
\left |\left\{0<\theta^{*}_{k}<1 \right\} \cap B_{3}^{*} \times [1,6]\right | > \delta,
\end{equation}
then 
\begin{equation}
\left |\left\{0<\theta^{*}<1 \right\} \cap B_{3}^{*} \times [1,6]\right | > (K_{1}+1)\delta \geq |B_{3}^{*} \times [1,6]|, 
\end{equation}
an absurdity.  Hence, there must exist a $k<K_{1}$ such that
\begin{equation}
\left |\left\{0<\theta^{*}_{k}<1 \right\} \cap B_{3}^{*} \times [1,6]\right | \leq \delta,
\end{equation}
and applying the Corollary \ref{Cor_Final_Isoperim} to $\theta_{k} \in \mathcal{S}(D,\alpha)$, we find that $\theta_{k}^{*} \leq 2-\lambda_{0}$ on $B_{1}^{*} \times [5,6]$, which translates to $\theta^{*} \leq 2-\lambda$ as desired.
\end{proof}
In the general case, we need to normalize and re-center $\theta$ in order to apply Claim \ref{eq:UnconditionDecreaseofLinf}.  Given an arbitrary $\theta \in \mathcal{S}(D,\alpha)$, define
\begin{align}
\overline{\theta}(t,x)=\frac{4}{\underset{[1,6] \times B_{3}^{*}}{\text{osc}} \theta^{*}} \left [ \theta(t,x)-\frac{1}{2} \left ( \inf_{[1,6] \times B_{3}^{*}} \theta^{*}+\sup_{[1,6] \times B_{3}^{*}} \theta^{*}\right )\right ].
\end{align}
We may assume 
\begin{equation}
\left |\left\{\overline{\theta}^{*}\leq 0\right \} \cap B_{3}^{*} \times [1,6] \right | \geq (1/2)\left | B_{3}^{*} \times [1,6]\right |.
\end{equation}
Indeed, if the opposite is true, we replace $\theta$ by $-\theta$ and the decrease in oscillations will come from below rather than above.  By design, $\overline{\theta}^{*}\leq 2$ on $B_{3}^{*} \times [1,6]$, so Claim \ref{eq:UnconditionDecreaseofLinf} implies that $\overline{\theta}^{*}\leq 2-\lambda$ in $B_{1}^{*} \times [5,6]$.  Translating this to $\theta^{*}$ gives 
\begin{align}
&\sup_{[5,6] \times B_{1}^{*}} \theta^{*} \leq \,\,\,\, (1/4)(2-\lambda)\underset{[1,6] \times B_{3}^{*}}{\text{osc}} \theta^{*} + \frac{1}{2} \left ( \inf_{[1,6] \times B_{3}^{*}} \theta^{*}+\sup_{[1,6] \times B_{3}^{*}} \theta^{*}\right). \\
\end{align}
Subtracting the infimum of $\theta^{*}$ over $[5,6] \times B_{1}^{*}$ from both sides and using elementary arguments, we may complete the proof with $\mu=1-\lambda/4$.
\end{proof}
We are now prepared to prove the main result, Theorem \ref{Thm:HolderContResult}.  In this final step, we need to repeatedly apply the oscillation reduction result \ref{cor:reductionofoscillation}, but on all scales.  It is here that we use the fact that the norm measuring the size of the velocity in the definition of $\mathcal{S}(D,\alpha)$ is critical.  This allows us to zoom in without changing $D$.
\begin{proof}[Proof of Main Theorem \ref{Thm:HolderContResult}]
Let $\theta \in \mathcal{S}(D,\alpha)$ be driven by a velocity field $u \in \LB{q}{-\gamma}$ where $2\alpha/q + \gamma = 2\alpha-1$.  

The first step of the proof is to apply Proposition \ref{cor:reductionofoscillation} on all scales and obtain a universal $C>0$ and $\beta \in (0,1)$ such that that for all $\theta \in \mathcal{S}(D,\alpha)$,
\begin{equation}\label{eq:ContinuityAtaPoint}
\sup_{(t,x) \in [5,6] \times B_{1}}\frac{|\theta(11/2,0)-\theta(t,x)|}{|t-11/2|^{\frac{\beta}{2\alpha}}+|x|^{\beta}}\leq C \|\theta\|_{L^{\infty}([1,\infty) \times \R^{d} )}.
\end{equation}
Towards this end, we define a scaling transformation
\begin{equation}
T(t,x,z)=[4^{-2}t + \frac{165}{32}, 4^{-\frac{1}{\alpha}}x,4^{-\frac{1}{\alpha}}z].
\end{equation}
The constants defining $T$ have been chosen to ensure two properties.  The first is that the vector $(11/2,0,0)$ is a fixed point of $T$.  The second is that $T(B_{3}^{*} \times [1,6]) \subset B_{1}^{*} \times [5,6]$.  In fact, for all $k \in \N$, it follows that $T^{k}(B_{3}^{*} \times [1,6]) \subset T^{k-1}\left (B_{1}^{*} \times [5,6] \right )$.  Since $\LB{q}{-\gamma}$ is a critical space, it follows that $\theta \circ T_{k} \in \mathcal{S}(D,\alpha)$ for all $k \in \N$.  Hence, applying the Proposition \ref{cor:reductionofoscillation} to $\theta \circ T_{k}$ gives
\begin{equation}
\underset{T_{k} \left ([5,6] \times B_{1}^{*} \right )}{\text{osc}}\theta^{*}\le \mu \underset{T_{k}([1,6] \times B_{3}^{*})}{\text{osc}} \theta^{*} \leq \mu \underset{T_{k-1}([5,6] \times B_{1}^{*})}{\text{osc}} \theta^{*}.
\end{equation}
Iterating this argument, we obtain the decay of oscillations
\begin{equation}
\underset{T_{k} \left ([5,6] \times B_{1}^{*} \right )}{\text{osc}}\theta^{*}\le \mu^{k} \underset{[1,6] \times B_{3}^{*}}{\text{osc}} \theta^{*}.
\end{equation}
Choosing $\beta$ such that $4^{\frac{\beta}{\alpha}}<1$, a small argument gives the estimate \eqref{eq:ContinuityAtaPoint}.  Since the H\"older estimate applies to all of $\mathcal{S}(D,\alpha)$ with a universal constant, re-centering in space and shifting forward in time gives
\begin{equation}
[\theta]_{C^{\beta}_{\alpha}([6,\infty) \times \R^{d})} \leq C \|\theta\|_{L^{\infty}([1,\infty) \times \R^{d} )},
\end{equation}
where the $C^{\beta}_{\alpha}$ semi-norm is defined by \eqref{eq:HolderSemiNorm}.  Finally, rescaling and using Lemma \ref{L.L2toLinfWholeSpace} gives
\begin{equation}
\|\theta\|_{C_{\alpha}^{\beta}([t,\infty) \times \R^{d})} \leqs t^{-(\beta + \frac{d}{4\alpha})}\|\theta_{0}\|_{L^{2}(\R^{d})},
\end{equation}
completing the proof.
\end{proof}
$\mathbf{Acknowledgements:}$
The authors wish to acknowledge that the Main theorem \ref{Thm:HolderContResult} was inspired by a conjecture in the list of open problems from the wiki page for Non-Local equations.  This work was partially supported by NSF grant DMS-1501067.

\bibliographystyle{plain}
\bibliography{Bibliography}

\end{document}